\documentclass[11pt,reqno]{article}
\usepackage{amssymb}
\usepackage{amsthm}
\usepackage{amsmath}
\newtheorem{theorem}{Theorem}[section]
\newtheorem{corollary}[theorem]{Corollary}

\newtheorem{example}[theorem]{Example}

\newtheorem{remark}[theorem]{Remark}

\begin{document}
\begin{center} \Large{\bf A point symmetry based method for transforming
ODEs with three-dimensional symmetry algebras  to their canonical forms}
\end{center}
\medskip

\begin{center}
 H. Azad$^*$, Ahmad Y. Al-Dweik$^*$, F. M. Mahomed$^{**}$ and M. T. Mustafa$^{***}$\\

{$^*$Department of Mathematics \& Statistics, King Fahd University
of Petroleum and Minerals, Dhahran 31261, Saudi Arabia}\\
{$^{**}$DST-NRF
Centre of Excellence in Mathematical and Statistical Sciences,
School of Computer Science and Applied Mathematics,
University of the Witwatersrand, Johannesburg, Wits 2050,  South Africa\\
{$^{***}$Department of Mathematics, Statistics and Physics, Qatar
University, Doha, 2713, State of Qatar}\\
}

hassanaz@kfupm.edu.sa, aydweik@kfupm.edu.sa, 
Fazal.Mahomed@wits.ac.za and tahir.mustafa@qu.edu.qa.
\end{center}
%%%%%%%%%%%%%%%%%%%%%%%%%%%%%%%%%%%%%%%%%%%%%%%%%%%%%%%%%%%%%%%%%%%%%%%%%%%%%%%%%%%%%%%%%%%%%%%%%%5
\begin{abstract}
We provide an algorithmic approach to the construction of point
transformations for scalar ordinary differential equations that
admit three-dimensional symmetry algebras which lead to their
respective canonical forms.
\end{abstract}
%%%%%%%%%%%%%%%%%%%%%%%%%%%%%%%%%%%%%%%%%%%%%%%%%%%%%%%%%%%%%%%%%%%%%%%%%%%%%%%%%%%%%%%%%%%%%
\section{Introduction}
This is a contribution to the algorithmic theory of
differential equations in the sense of Schwarz [6]. 
In the first five sections we provide results needed to construct algorithms, which are illustrated 
in detail in section 6 of the paper. 

Not all the three-dimensional algebras have invariant differential equations of
second-order in the sense that they are the full symmetry algebra
of the  equation. For this reason, we have given  higher-order
invariant equations for such types of three-dimensional algebras.

Although all the results of this paper were independently arrived at, 
the main ideas are already in Lie \cite{Lie1891}. The justification for 
placing this in the public domain is its brevity, clarity and a uniform treatment 
of compact and non-compact algebras.

The realizations of two- and three-dimensional Lie algebras as
vector fields in $\mathbb{C}^2$ is given in e.g., Ibragimov
\cite[page 163]{Ibragimov1999}. The details, as well as the
invariant second-order ordinary differential equations (ODEs) are
given in Lie \cite[pages 479-542]{Lie1891}.
The two-dimensional algebras are essentially distinguished by
their ranks. It is thus desirable to give a similar description of
three-dimensional Lie algebras over the reals.

The main aim of this work is to give such a description and an
algorithmic procedure that systematically utilizes the structural
information more explicitly and which is programmable. This
program reduces any given ODE that admits a three-dimensional
algebra to its simplest form.

This is achieved by proving a version of the Lie-Bianchi
classification in an algorithmic way and giving the realizations
of the algebras as vector fields in $\mathbb{R}^2$.

In every case, there is an invariantly defined two-dimensional
abelian algebra - which is not always a subalgebra of the given
algebra - and its rank determines the coordinates that reduces the
equation to its simplest form.

As far as the form of the invariant differential equations is
concerned, one can use the result in \cite[pages 69-76]{Azad2015}
that reduces the computation of local joint invariants of any
finite number of vector fields algorithmically to that of an
abelian algebra of vector fields.

A few words regarding the importance of low-dimensional algebras
for differential equations are in order.

Lie [3,4] obtained the complete explicit classification of scalar
second-order ODEs that possess
non-similar (not transformable into each other via point
transformation) complex Lie algebras of dimension $r$, where
$r=0,\ldots,8$. He showed that the complex Lie algebra of vector
fields acting in the plane admitted by a given second-order ODE
can only be of dimensions $0,1,2,3$ or $8$. He also proved that if
a second-order equation admits an eight-dimensional algebra, it is
linearizable by means of a point transformation and it is then
equivalent to the simplest equation, viz. the free particle equation.

It is well-known that one- and two-dimensional algebras have
identical structures over the reals as well as over the complex
numbers. As a consequence, the Lie symmetry algebra classification
of scalar second-order ODEs over the reals is precisely the same
as that over the complex numbers for one- and two-dimensional Lie-algebras.

If a second-order equation admits a single generator of symmetry,
then in general its order can be reduced by one [4]. Moreover, Lie
[4] showed that scalar second-order equations possessing two
generators of symmetry have four canonical forms. These are
well-known now as the Lie canonical forms for the vector fields
and their representative second-order equations. Lie [4] also
proved that the rank one algebras result in linearization of the
associated second-order ODE.

The situation is different for three- and higher-dimensional Lie
algebras as there are fewer complex than real algebras of
dimension three or more. This arose in the Bianchi [2]
classification of Lie algebras. Two of the complex Lie algebras of
dimension three in the real plane each split up into two real
non-isomorphic  Lie algebras. Therefore, there are two more
non-isomorphic real three-dimensional Lie algebras than complex algebras.

Due to the above considerations on Lie algebras in the real plane
for higher dimensions, there are additional three-dimensional
algebras of vector fields acting in the real plane than in the
complex plane. These were deduced by Mahomed and Leach [5]. These
yield additional non-similar scalar second-order equations that
admit real Lie algebras [5].

In summary, this is a contribution to the algorithmic Lie theory
of scalar ODEs in the sense of Schwarz
[6]. Schwarz [6] utilized janet bases in the representation of the
determining equations of the symmetry generators. The main
difference in our approach is to use canonical forms of the
symmetry algebra in order to construct the requisite point
transformations that bring a given ODE
with known symmetry algebra to its canonical form.

The reader is referred to Ibragimov [7] for the background on
Lie's theory of symmetries of differential equations.
%%%%%%%%%%%%%%%%%%%%%%%%%%%%%%%%%%%%%%%%%%%%%%%%%%%%%%%%%%%%%%%%%%%%%%%%%%%%%%%%%%%%%%%%%%%%%%
\section{The Lie-Bianchi Classification of three-dimensional solvable algebras}

We begin with a formulation of the Lie-Bianchi [1,2]
classification of three-dimensional Lie algebras.

\begin{theorem}\label{th111}
 Let $G$ be a three-dimensional Lie algebra. If $G$ is
solvable, then $G'$ is abelian. Moreover,

(a) if $G'$ is two-dimensional, then the structure of $G$ is
completely determined by the eigenvalues and multiplicities of
{\rm ad} $(X)$ as a linear transformation of $G'$, where $X$ is any
representative of
$G/G'$ in $G$,

(b) if $G'$ is one-dimensional, then the structure of $G$ is
completely determined by the dimension of the centralizer of $G'$
in $G$.
\end{theorem}

\proof In general, by Lie's theorem on complex solvable linear Lie
algebras \cite[page 106]{Neeb2012}, if $G$ is a solvable Lie
algebra then the algebra ad $(G^{\mathbb{C}})$ is nilpotent, where
$G^{\mathbb{C}}=G+\sqrt{-1}~G$. Therefore, the commutator $G'$ of
$G$ is nilpotent.

Now assume that $G$ is solvable and of dimension three.

(a) If $H$ is a two-dimensional algebra with basis $\{X,Y\}$, then its
commutator is generated by $[X,Y]$. If $H$ is nonabelian, extend
$U=[X,Y]$ to a basis $\{U,V\}$ of $H$. Then, scaling $V$, we have
the canonical representation of $H$ by the relations $[V,U]=U$; such
an algebra is not nilpotent.

Therefore if $G'$ is two-dimensional, it must be abelian. Take a
basis $\{X,Y\}$ of $G'$ and extend it to a basis $\{X,Y,Z\}$ of $G'$. Then
ad $(Z)$ operating on $G'$ does not have 0 as an eigenvalue and
the eigenvalues and their multiplicities of ad $(Z)$ operating on
$G'$ completely determine the structure of $G$.

(b) Assume that $G'$ is one-dimensional. Let $\{U\}$ be a basis of $G'$.
Extend it to a basis of $\{X,Y,U\}$ of $G$.

So $[X,Y]=aU, [X,U]=bU,[Y,U]=cU$. Now dim $Z_{A}(U)\geq 1.$
Suppose it is one. Then $b,c \neq 0$. By scaling $X,Y$ suitably, we
then have $[X,Y]=aU,[X,U]=U,[Y,U]=U$ so $[X-Y,U]=0$. Therefore
$X-Y, U$ are in $Z_G(U)$ and dim $Z_G (U)$ is at least two.

(i) Suppose dim $Z_{G}(U)=2.$ Choose a basis $\{Y,U\}$ of $Z_{G}(U)$
and extend it to a basis $\{X,Y,U\}$ of $G$. So
$[X,Y]=aU,[X,U]=bU,[Y,U]=0$, with $b\neq 0$. By scaling, we have
that $[X,Y]=aU,[X,U]=U,[Y,U]=0$.

Now, for any $\lambda$ we have
\[[X,Y+\lambda U]=aU+\lambda U, [X,U]=U, [Y+\lambda U,U]=0.\]
Choosing $a+\lambda =0$ and renaming $Y-aU$ as $Y$ we have the
canonical relations $[X,Y]=0,[X,U]=U,[Y,U]=0$ and
$Z_{G}(U)=\langle U,Y\rangle$.

(ii) Suppose dim $Z_{G}(U)=3$. Take any basis, $\{X,Y,U\}$ of $G$. We
then have
\[[X,Y]=aU,[X,U]=0,[Y,U]=0, \,\,\mbox{with}\, a\neq 0.\]
By scaling, we can assume that $[X,Y]=U,[X,U]=0,[Y,U]=0$.

Therefore the canonical representations of $G$ when dim $G'$ is one
are:
\[ [X,Y]=0,[X,U]=U,[Y,U]=0\,\,\mbox{(dim}\,Z_{G}(G')\,\mbox{is}\,
2)\]
and 
\[ [X,Y]=U,[X,U]=0,[Y,U]=0\,\,\mbox{(dim}\,
Z_{G}(G')\,\mbox{is}\, 3).\]
\endproof

\begin{corollary} If $G'$ is two-dimensional then it is abelian
and there is a basis $\{X,Y,Z\}$ of $G$ with $\{X,Y\}$ a basis of $G'$
with {\rm ad} $(Z)$ operating on $G'$ as follows - recall that {\rm ad} $(Z)$ does
not have $0$ as an eigenvalue of $G'$:

(i) {\rm ad} $(Z)$ has real and distinct eigenvalues.

If $X,Y$ are eigenvectors, then as {\rm ad} ($Z$) does not have 0 as an
eigenvalue of $G'$, by scaling $Z$ we have
$[Z,X]=X,[Z,Y]=cY,c\neq 0,1$.

(ii) {\rm ad} $(Z)$ has only one eigenvalue and the corresponding
eigenspace is two-dimensional.

In this case, scaling $Z$, we have the canonical representation
\[[Z,X]=X,[Z,Y]=Y.\]

(iii) {\rm ad} $(Z)$ has only one real eigenvalue and the corresponding
eigenspace is one-dimensional. There is a basis $\{X,Y\}$ of $G'$ with
$[Z,X]=X+Y,[Z,Y]=Y$.

(iv) {\rm ad} $(Z)$  has a complex eigenvalue $\lambda$. The canonical
relations - after scaling $Z$ - are
\[[Z,X]=\cos \theta ~X+ \sin \theta ~Y,[Z,Y]=-\sin \theta ~X+\cos \theta ~Y,\]
where $X,Y$ are the real and imaginary parts of an eigenvector for the
eigenvalue $\lambda.$
\end{corollary}

\proof Only cases (iii) and (iv) need  proofs.

Case (iii): ad $(Z)$ has only one real eigenvalue, which is non-zero
and the null-space of ad $(Z)-\lambda I$ is one-dimensional. Let
$V$ be an eigenvector for the eigenvalue $\lambda$. Extend $V$ to a
basis $\{X,V\}$ of $G'$. Then $X$ is a generalized eigenvector of $Z$
and $\{X,Y=$ ad $(Z-\lambda I)\,(X)\}$ is a basis of $G'$.

We have $[Z,X]=\lambda X+Y,[Z,Y]=\lambda Y$. Dividing $Z$ by
$\lambda$ and relabeling it $Z$ we have the relations
$[Z,X]=X+cY,(c\neq 0),[Z,Y]=Y$. Finally, replacing $Y$ by $cY$, we
have the canonical relations $[Z,X]=X+Y,[Z,Y]=Y$.

Case (iv): if one of the eigenvalues of ad ($Z$) is complex but not
real, then the other eigenvalue is its conjugate. Denote by
$\lambda$ any one of these eigenvalues. Then the eigenvectors live
in the Lie algebra $G+\sqrt{-1}~G$ - the complexification of $G$.
We will denote the complexification of any Lie algebra $G$ by
$G^{\mathbb{C}}$. The eigenvectors of ad $(Z)$ live in
$G^{'\mathbb{C}}$. Find an eigenvector $e$ for ad $(Z$).

The real and imaginary parts of $e$ are $\displaystyle
{\rm Re}(e)=\frac{e+\bar{e}}{2}, {\rm Im}(e)= \frac{e-\bar{e}}{2i}$.

By scaling $Z$ by $\frac{1}{|\lambda|}$, we have the canonical
relations
\[[Z, {\rm Re}(e)]=\cos \theta \,{\rm Re}(e)+\sin \theta \,{\rm Im}(e)\]
\[[Z, {\rm Im}(e)]=-\sin \theta \, {\rm Re}(e)+\cos \theta \, {\rm Im}(e).\]
Take $X={\rm Re}(e)$ and $Y={\rm Im}(e)$.
\endproof
%%%%%%%%%%%%%%%%%%%%%%%%%%%%%%%%%%%%%%%%%%%%%%%%%%%%%%%%%%%%%%%%%%%%%%%%%%%%%%%%%%%%%%%%%%%%%%%%%%%%%5
\section{Local classification of  commuting Lie algebras of vector
fields in $\mathbb{R}^2$}
 
The algorithms for bringing a given ODE
to its canonical form by point transformations ultimately reduce
to constructive classifications of abelian Lie algebras of vectors
fields. This section is devoted to an algorithmic procedure for
finding canonical forms of such algebras.

It is well-known that  if $X$ is a vector field and $X(p)\neq 0$,
then  near $p$ we can introduce coordinates in which
$X=\partial_{x}$. To find such a canonical coordinate, one uses the
method of characteristics \cite[page 142]{Ibragimov1999}, to find
a basic invariant function $y$ of $X$ and choose a function $x$
functionally independent from $y$. Then in these coordinates
$X=f(x,y)\,\partial_{x},$ with $f$ non-vanishing near $p$.

We want functions $\tilde{x},\tilde{y}$ with $X(\tilde{x})=1,\,
X(\tilde{y})=0$. Then the requirements become
$f(x,y)\,\partial_{x}(\tilde{x})=1,
f(x,y)\,\partial_{x}\,(\tilde{y})=0.$
A solution of this system is $\displaystyle \tilde{x}=\int
\frac{dx}{f(x,y)},\,\tilde{y}=y$.

Now suppose that $X,Y$ are commuting vector fields with rank 2
near a point $p$; so that $X(q),Y(q)$ are linearly independent near
$p$; say $X(p)\neq 0$. We may therefore assume that in some
neighbourhood of
$p, X(q)\neq 0$ and $X,Y$ of rank 2 in this neighbourhood.

Choose local coordinates $x,y$ in  possibly a smaller neighbourhood
with $X=\partial_{x}$. As $Y$ operates on invariants of
$X$, $Y(y)=g(y)$. Now as $Y$ operates non-trivially on invariants of $X$
because of the rank condition, we can find a function $\tilde{y}$
of $y$ with $Y(\tilde{y})=1$. By change of notation, we now have
coordinates $x,y$ with $X=\partial_{x},Y=\xi
(y)\,\partial_{x}+\partial_{y}$. 
We want new coordinates $\tilde{x},\tilde{y}$ with
$X(\tilde{x})=1,
X(\tilde{y})=0, \, Y(\tilde{x})=0,\, Y(\tilde{y})=1.$
The system to solve now is $\displaystyle \frac{\partial
\tilde{x}}{\partial x}=1,\frac{\partial \tilde{y}}{\partial x}
=0\,\,\xi (y)+\frac{\partial \tilde{x}}{\partial
y}=0,\frac{\partial
\tilde{y}}{\partial y}=1.$
The solution is given by $\displaystyle \tilde{x}=x+\varphi (y),
\, \frac{\partial \tilde{x}}{\partial y}=\varphi'(y)=-\xi
(y),\,\tilde{y}=y.$
If the fields are commuting of rank 1, and $X=\partial_{x}$ in
local coordinates, then $Y=f(y)X,$ with $f'(y)\neq 0$. Then in the
variables $x,\tilde{y}=f(y)$ the canonical form of the fields is
$X=\partial_{x},
Y=\tilde{y}\,\partial_{x}$.
%%%%%%%%%%%%%%%%%%%%%%%%%%%%%%%%%%%%%%%%%%%%%%%%%%%%%%%%%%%%%%%%%%%%%%%%%%%%%%%%%%%%%%%%%%%%%%%%%%%%%%%%%%%%%%%%%%%%%%%
\section{ Realizations of three-dimensional  algebras as vector fields in $\mathbb{R}^2$ }

We now turn to realizations of the algebras occurring in Theorem
\ref{th111} and its corollary as algebras of vector fields in
$\mathbb{R}^2$. This was done by Lie for complex Lie algebras that
arise as symmetries of second-order ODEs: see \cite[page
164]{Ibragimov1999} and \cite[pages 479-530]{Lie1891}. Mahomed and
Leach [5] extended this study for real Lie
algebras. Similar ideas apply to vector fields in $\mathbb{R}^3$.

$\mathbf{Outline~of~the~argument}$: Using Theorem \ref{th111}, the
canonical forms for solvable three-dimensional algebras with
nontrivial commutator can be easily obtained algorithmically
as one has just to put a
two-dimensional algebra, depending on its rank, in canonical form.
The reason is that if $G'$ is two-dimensional then it is abelian and
its rank determines the canonical form of $G$. If $G'$ is of
dimension one and its centralizer has dimension two, then the rank of
the centralizer determines the canonical form of $G$. Finally, if
$G'$ is one-dimensional and its centralizer has dimension three, then
rank of $G$ must be 2 and choosing a field supported outside of
$G'$ gives an abelian two-dimensional algebra of rank 2 which determines the
canonical form of $G$.

If $G'=G$, then picking any nonzero element $X$ of $G$, the
eigenvalues of ad $(X$) determine the canonical form of $G$ and of
the vector fields also - as detailed below in section \ref{4.3}.

$\mathbf{Details~of~the~classification}$: Assume that $G$ is 
three-dimensional and abelian. Then its rank is at most 2 and it cannot
be 2 as the centralizer of $\{
\partial_{x},\partial_{y}\}$ is $\langle
\partial_{x},\partial_{y}\rangle$. Therefore it is of rank 1.
Pick any non-zero element $X$ of $G$ and find canonical
coordinates
$x,y$ so that $X=\partial_{x}$ and extend it to a basis $\{X,Y,Z\}$.
Since rank of $G$ is 1, necessarily
$Y=f(y)\,\partial_{x},\,Z=g(y)\,\partial_{x}$.
Since $f$ is not a constant, we can take $\tilde{y}=f(y)$.
So the basis becomes
$\{\partial_{x},\tilde{y}\partial_{x},\,g(f^{-1}\,
(\tilde{y}))\,\partial_{x}=h(\tilde{y})\,\partial_{x}\}$,
where $h(\tilde{y})$ is linearly independent of $\{1,\tilde{y}\}$.

For the canonical realizations of the algebras occurring in
Theorem 1, one needs to solve equations of the type $[Z,X]=aX+bY,
\,[Z,Y]=cX+dY$, where $X=\partial_{x}$ and $Y=\partial_{y}$ or $y
\,\partial_{x}$.
In case $Y=\partial_{y}$,  this is  straightforward. However, when
$Y=y\partial_{x}$, we have $Z=\xi \partial_{x}+\eta
\,\partial_{y}$, where
\[\xi =-ax-bxy+\varphi (y), \,\eta =c+(d-a)y-by^2.\]
We want a change of variables $\tilde{x},\tilde{y}$ so that
$\partial_{x}=\partial_{\tilde{x}},\,
y\,\partial_{x}=\tilde{y}\,\partial_{\tilde{x}}$ and
\begin{equation}
Z=(-a\tilde{x}-b\tilde{x}\,\tilde{y})\,\partial_{\tilde{x}}+(c+(d-a)\,\tilde{y}-b\tilde{y}^{2})\, \partial_{\tilde{y}}\label{*}
\end{equation}
Then necessarily $\tilde{y}=y,\,\tilde{x}=x+\psi (y),
\,\partial_{\tilde{y}}=-\psi^{'}(y)\,\partial_{x}+\partial_{y}$.
Substituting these expressions in the equation (\ref{*}) we arrive at
equation
\begin{equation}
\varphi (y)=\psi(y)(-a-by)-\psi^{'}(y)(c+(d-a)y-by^2)\label{**}
\end{equation}
Solving this differential equation for $\psi$ removes the term
$\varphi (y)$ in the field $Z$ in the new variables.

This gives the following realizations of the algebras occurring in
Theorem 1 as vector fields in the plane:
%%%%%%%%%%%%%%%%%%%%%%%%%%%%%%%%%%%%%%%%%%%%%%%%%%%%%%%%%%%%%%%%%%%%%%%%%%%%%%%%%%%%%%%%%%%%%%
\subsection{${\rm dim}\,(G\,')=1$}
I: $\mbox{dim}\, Z_G(G\,')=2,\,[X,Y]=0,\,[X,U]=U,[Y,U]=0,$ rank
$Z_{G}\,(G\,')=1$\\
$U=\partial_{y}, Y=f(x)\,\partial_{y},$ where $f$ is not a
constant. Making a change of variables $\tilde{x}=f(x),$
relabeling $\tilde{x}=x,$ using (\ref{*}), (\ref{**})- with $Z$ replaced by
$X$ - we obtain the representation
\[U=\partial_{y}, Y=x\,\partial_{y}, X=-y\partial_{y}-x\,\partial_{x}\]

II: $\mbox{dim}\, Z_G(G\,')=2,\,[X,Y]=0,\,[X,U]=U,[Y,U]=0,$ rank
$Z_{G}\,(G\,')=2,$
\[U=\partial_{y},\, Y=\partial_{x}, X=-y\partial_{y}\]

III: dim $Z_G(G\,')=3,[X,Y]=U,[X,U]=0,[Y,U]=0$ \\
In this case rank of $G$ must be 2.
In the canonical coordinates for $U$, we have $U=\partial_{y},$
and one of $X$ or $Y$ is supported outside $\partial_{y},$
otherwise $[X,Y]$ would be 0. By symmetry between $X$ and $Y$, we
may suppose that $X$ is supported outside $\partial_{y}$ and
therefore rank
$\langle X,U\rangle=2.$
We may then suppose that $X=\partial_{x}.$ This determines
$Y=x\,\partial_{y}.$ The canonical realization is thus
$X=\partial_{x}, U=\partial_{y},\, Y=x\,\partial_{y}.$
%%%%%%%%%%%%%%%%%%%%%%%%%%%%%%%%%%%%%%%%%%%%%%%%%%%%%%%%%%%%%%%%%%%%%%%%%%%%%%%%%%%%%%%%%%%%%%
\subsection{${\rm dim}\,(G\,')=2$}
(a) There is a basis $\{X,Y,Z\}$ of $G$ with $\{X,Y\}$ a basis of $G'$
with ad $(Z)$ operating on $G\,'$ having real and distinct eigenvalues.
Computing the eigenvectors for these eigenvalues and labeling
them $X,Y$ and dividing $Z$ by the eigenvalue for the eigenvector
$X$, and relabeling it $Z$ we have the relations
$[Z,X]=X, [Z,Y]=cY,\, c\neq 0,1.$

IV: rank $(G\,')=2.$\\
 Choosing coordinates with $X=\partial_{x},\,
Y=\partial_{y}$ we get $Z=-x\,\partial_{x}-cy\,\partial_{y}$

V: rank $(G\,')=1.$\\
 Choosing coordinates in which
$X=\partial_{y},\, Y=x\,\partial_{y}$ and using (\ref{*}) and (\ref{**}) we
obtain - by change of
notation $Z=(c-1)\, x\,\partial_{x}-y\,\partial_{y}$.

(b) There is a basis $\{X,Y,Z\}$ of $G$ with $\{X,Y\}$ a basis of $G\,'$
with ad $(Z)$ operating on $G\,'$ having a real eigenvalue with the
corresponding eigenspace of dimension two.
Computing the eigenvectors for these eigenvalues and labeling
them $X,Y$ and dividing $Z$ by the eigenvalue and relabeling it
$Z$ we have the relations
\[[Z,X]=X,\,[Z,Y]=Y.\]

VI: rank $(G\,')=2.$\\
 Choosing coordinates with
$X=\partial_{x},\,Y=\partial_{y}$ we deduce
$Z=-x\,\partial_{x}-y\,\partial y$

VII: rank $(G\,')=1.$\\ 
Choosing coordinates with
$X=\partial_{y},\,Y=x\,\partial_{y}$ and using (\ref{*}) and (\ref{**}) we get - by change of notation - $Z=-y\,\partial_{y}.$

(c) There is a basis $\{X,Y,Z\}$ of $G$ with $\{X,Y\}$ a basis of $G\,'$
with ad $(Z)$ operating on $G\,'$ having only one eigenvalue $\lambda$
with the corresponding eigenspace of dimension one.
Find the corresponding eigenvector and let $X$ be a vector
linearly independent from this eigenvector. Then $X$ is a
generalized eigenvector and $\{X,Y=(\mbox{ad}\,(Z)-\lambda I)\,X\}$ is a
basis of $G'$. Dividing $Z$ by $\lambda$ and relabeling it $Z$, we
have
$ [Z,X]=X+\frac{1}{\lambda}\, Y, [Z,Y]=Y.$
Replacing $Y$ by $\frac{1}{\lambda}Y$ in the second
equation and finally denoting $ \frac{1}{\lambda}Y$ by $Y$ we have the relations $[Z,X]=X+Y,\,[Z,Y]=Y,\,[X,Y]=0.$

VIII: rank $(G\,')=2.$\\
 Choosing coordinates with
$X=\partial_{x},\, Y=\partial_{y}$ we arrive at
\[Z=-x\,\partial_{x}-(x+y)\,\partial_{y}\]

IX: rank $(G\,')=1.$\\
 Choosing coordinates in which
$X=x\,\partial_{y},\, Y=\partial_{y}$ and using (\ref{*}) and (\ref{**}) we
find - by change of notation - $Z=\partial_{x}-y\,\partial_{y}$.

(d) There is a basis $\{X,Y,Z\}$ of $G$ with $\{X,Y\}$ a basis of $G\,'$
and ad  $(Z)$ operating on $G\,'$ having a
non-real complex eigenvalue $\lambda.$ Let $e$ be an eigenvector.
Dividing $Z$ by $\frac{1}{|\lambda|},$ denoting it again by $Z$ we
have the canonical relations
\[[Z,\,{\rm Re}(e)]=\cos \theta \,{\rm Re}(e)-\sin \theta \, {\rm Im}(e)\]
\[[Z,{\rm Im}(e)]=\sin \theta \,{\rm Re}(e)+\cos \theta \, {\rm Im}(e).\]
Take $X={\rm Re}(e)$ and $Y={\rm Im}(e).$

X: rank $(G\,')=2.$\\
 Choosing coordinates with $X=\partial_{x}, \,
Y=\partial_{y}.$ Then
\[Z=(-x \cos \theta - y \sin \theta)\,\partial_{x}+(x\,\sin \theta - y \, \cos \theta)\,\partial_{y}.\]

XI: rank $(G\,')=1.$\\
Choosing coordinates in which
$X=\partial_{y},\,Y=x\,\partial_{y}$ and using (\ref{*}) and (\ref{**}) we
obtain - by change of notation - $Z=(1+x^2)\,\sin \theta \,
\partial_{x}+y( x\sin \theta-\cos\theta )\,\partial_{y}$
%%%%%%%%%%%%%%%%%%%%%%%%%%%%%%%%%%%%%%%%%%%%%%%%%%%%%%%%%%%%%%%%%%%%%%%%%%%%%%%%%%%%%%%%%%%%%%
\subsection{${\rm dim}\, (G\,')=3$}\label{4.3}
Let $G$ be a three-dimensional Lie algebra with $G'=G$. If $I$ is
any non-zero vector in $G$, then its centralizer consists of
multiples of $I$; for if $U$ is linearly independent of $I$ and it
centralizes $I$, then extending $I$ to a basis of $G$ we see that
$G'$ is at most two-dimensional. Therefore, the non-zero
eigenvalues of ad $(I)$ occur in pairs $\lambda,-\lambda$ - as the
trace of ad $(I)$ is zero.

Case (i): The non-zero eigenvalues of ad $(I)$ are real given as $\pm \lambda$.
Then we can find eigenvectors $U,V$ of $I$  with $[I,U]=\lambda~
U,[I,V]=-\lambda V,[U,V]=cI$ where $c\ne 0$ as $[U,V]$
centralizes $I$. Setting
$X=U,~Y=\frac{2}{c~\lambda}~V,~Z=\frac{2}{\lambda}~I$, we have the
standard relations of $sl(2,\mathbb{R})$ given by
$$[Z,X]=2X,[Z,Y]=-2Y,[X,Y]=Z.$$
In this case the Killing form is non-degenerate and indefinite.

Case (ii): ad $(I)$ has a non-real eigenvalue.
In this case, the nonreal eigenvalues must be purely imaginary.
Scaling $I$, we may suppose that the non-zero eigenvalues of
ad $(I)$ are $\pm \sqrt{-1}$. Take an eigenvector $e$ of ad $(I)$ in
the complexification of $G$ with eigenvalue $\sqrt{-1}$. Write
$e=U-\sqrt{-1}~V$. Then $[I,U]=V$ and $[I,V]=-U.$ Now
$[e,\bar{e}]$ commutes with $I$ and
$[e,\bar{e}]=2~\sqrt{-1}~[U,V]$. Therefore $[U,V]$ commutes with
$I$. We can scale these generators such that $[I,U]=V,[I,V]=-U$
and $[U,V]=\epsilon I$, where $\epsilon^2=1$. If $\epsilon=-1$,
then we have the generators with $[I,U]=V,[I,V]=-U$ and
$[U,V]=-I$. So $[U,V+I]=-(V+I)$ and $[U,V-I]=V-I$ and we are back
to case (i). In this case the
Killing form is non-degenerate and indefinite. 
If $\epsilon=1$, then  we have the generators with
$[I,U]=V,[I,V]=-U$ and $[U,V]=I$. The Killing form is  negative definite. These are the standard relations of $so(3)$.

Case (iii): All the eigenvalues of ad $(I)$ are zero.
In this case ad $(I)$ is nilpotent and in the normalizer $N(I)$ of
$I$ there must be an element with real non-zero eigenvalues,
specifically, any element $H$ in $N(I)$ complementary to $I$ must
have non-zero eigenvalues, so $[H,I]=\lambda I,[H,U]=-\lambda U$
for some element $U$ and we are back to case (i).
\begin{remark}
If the Killing form is definite, then the Lie algebra must be
$so(3)$. So its Killing form must in any case be negative
definite. In this case the non-zero eigenvalues of {\rm ad} $(I)$ must be
purely imaginary. Arguing exactly as above, given a non-zero
element $I$ of $L$, we can find generators $U,V$ with
$[I,U]=V,[I,V]=-U$ and $[U,V]=cI$, where $c=\pm 1$. Now $c=-1$
would give an indefinite Killing form, so $[I,U]=V,[I,V]=-U$ and
$[U,V]=I$ which give the relations for $so(3)$ for the triple
$U,V,I$. In this case clearly
$[I,\sqrt{-1}~U]=\sqrt{-1}~V,[I,\sqrt{-1}~V]=-\sqrt{-1}~U$ and
$[\sqrt{-1}~U,\sqrt{-1}~V]=-I$ which are the relations for
$sl(2,\mathbb{C})$. On the other hand, if $[I,U]=
U,[I,V]=-V,[U,V]=I,$ the non-zero eigenvalues of {\rm ad} $(U+V)$ are
purely imaginary and working with the corresponding eigenvectors we
obtain a basis for $so(3,\mathbb{C})$. For this reason, over the
complex numbers, there is a correspondence between
$sl(2,\mathbb{R})$ and $so(3)$ invariant equations.
\end{remark}
%%%%%%%%%%%%%%%%%%%%%%%%%%%%%%%%%%%%%%%%%%%%%%%%%%%%%%%%%%%%%%%%%%%%%%%%%%%%%%%%%%%%%%%%%%%%%%
\subsubsection{Realizations of $sl(2,\mathbb{R})$ as vector fields in ${\mathbb{R}}^2$} 
Here we discuss the indefinite case.
So let $L$ be a three-dimensional algebra of vector fields with
indefinite Killing form. Then $L$ has a basis $\{X, Y, Z\}$ such that
$$[Z,X]=2X,[Z,Y]=-2Y,[X,Y]=Z.$$
Find coordinates in which $Z=\partial_x$. Then necessarily $X$ has
the following form
\begin{equation}\label{te6}
\begin{array}{cc}
X=e^{2 x}\left(  f_1(y)\partial_{x}+f_2(y)\partial_{y} \right).\\
\end{array}
\end{equation}
As $\{\partial_x, f_1(y)\partial_{x}+f_2(y)\partial_{y}\}$ are
commuting vectors, we have the following two cases:

Case A: the rank of $\{\partial_x, f_1(y)\partial_{x}+f_2(y)\partial_{y}\}$ is 2.
There is a change of variables in which
$\partial_x=\partial_{\tilde{x}}$ and
$f_1(y)\partial_{x}+f_2(y)\partial_{y}=\partial_{\tilde{y}}$.
Using section 2, such a change of variables can be given explicitly
as follows:
\begin{equation}
\begin{array}{lll}
\tilde{x}=x-\int{\frac{f_1}{f_2}dy},& \tilde{y}&=\int{\frac{1}{f_2}dy}.\\
\end{array}
\end{equation}
Therefore in these coordinates
\begin{equation}\label{te7}
\begin{array}{cc}
X=e^{2 (\tilde{x}+f(\tilde{y}))} \partial_{\tilde{y}},\\
\end{array}
\end{equation}
where $f(\tilde{y})=\int{\frac{f_1}{f_2}dy}.$ In the coordinates
$\tilde{\tilde{x}}$, $\tilde{\tilde{y}}$ in which
$\tilde{\tilde{x}}=\tilde{x}$ and $e^{2 f(\tilde{y})}
\partial_{\tilde{y}}=\partial_{\tilde{\tilde{y}}}$ given by
\begin{equation}
\begin{array}{lll}
\tilde{\tilde{x}}=\tilde{x},& \tilde{\tilde{y}}&=\int{e^{-2 f(\tilde{y})} d\tilde{y}},\\
\end{array}
\end{equation}
we have
\begin{equation}\label{te8}
\begin{array}{cccc}
Z=\partial_{\tilde{\tilde{x}}},&X=e^{2 \tilde{\tilde{x}}} \partial_{\tilde{\tilde{y}}},&Y=e^{-2\tilde{\tilde{ x}}}\left(  (\tilde{\tilde{y}}+c_1)\partial_{\tilde{\tilde{x}}}+\left( {(\tilde{\tilde{y}}+c_1)}^2+\epsilon~\lambda^2\right)\partial_{\tilde{\tilde{y}}} \right),\\
\end{array}
\end{equation}
where $\epsilon \in \{0,1,-1\}$.
So, the new coordinates
$\bar{x}=\tilde{\tilde{x}},~\bar{y}=\frac{1}{\lambda}(\tilde{\tilde{y}}+c_1)$
transform $X,Y$ and $Z$ to
\begin{equation}\label{te9}
\begin{array}{cccc}
Z=\partial_{\bar{x}},&X=\frac{1}{\lambda}e^{2 \bar{x}} \partial_{\bar{y}},&Y=\lambda e^{-2\bar{x}}\left(  \bar{y}\partial_{\bar{x}}+\left( {\bar{y}}^2+\epsilon\right)\partial_{\bar{y}} \right),\\
\end{array}
\end{equation}
where $\epsilon \in \{0,1,-1\}$.
Finally, the transformation $\hat{x}=e^{-2\bar{x}},\hat{y}=
\bar{y} e^{-2\bar{x}}$ transforms $X,Y$ and $Z$ to the polynomial
form
$$ X=\partial_{\hat{y}}, Y=-2\hat{x}\hat{y}\partial_{\hat{x}}+({-\hat{y}}^2+\epsilon~{\hat{x}}^2)\partial_{\hat{y}}, Z= -2\hat{x}\partial_{\hat{x}}-2\hat{y}\partial_{\hat{y}},$$
where $\epsilon \in \{0,1,-1\}$.

Case B: the rank of $\{\partial_x, f_1(y)\partial_{x}+f_2(y)\partial_{y}\}$ is 1.
In this case the vector $X$ should have the form
\begin{equation}\label{te10}
\begin{array}{cc}
X=e^{2 x} f_1(y)\partial_{x}.\\
\end{array}
\end{equation}
If $f_1(y)$ is constant, then
\begin{equation}
\begin{array}{cccc}
Z=\partial_{x},&X=e^{2 x} \partial_{x}.\\
\end{array}
\end{equation}
Using $[X,Y]=Z$, gives
\begin{equation}
\begin{array}{cccc}
Y=-\frac{1}{4}e^{-2x} \partial_{x}.\\
\end{array}
\end{equation}
If $f_1(y)$ is not a constant, we can introduce a change of
variables
\begin{equation}
\begin{array}{lll}
\tilde{x}=x,& \tilde{y}&=f_1(y).\\
\end{array}
\end{equation}
Therefore in these coordinates
\begin{equation}\label{te11}
\begin{array}{cc}
X=e^{2 (\tilde{x}+f(\tilde{y}))} \partial_{\tilde{x}},\\
\end{array}
\end{equation}
where $f(\tilde{y})=\frac{1}{2}\ln{\tilde{y}}.$
Finally, in the coordinates $\tilde{\tilde{x}}$,
$\tilde{\tilde{y}}$ given by
\begin{equation}
\begin{array}{lll}
\tilde{\tilde{x}}=\tilde{x}+\frac{1}{2}\ln{\tilde{y}},& \tilde{\tilde{y}}&=\tilde{y},\\
\end{array}
\end{equation}
we have
\begin{equation}\label{te12}
\begin{array}{cccc}
Z=\partial_{\tilde{\tilde{x}}},&X=e^{2 \tilde{\tilde{x}}} \partial_{\tilde{\tilde{x}}}.\\
\end{array}
\end{equation}
Using $[X,Y]=Z$, gives
\begin{equation}\label{te13}
\begin{array}{cccc}
Y=-\frac{1}{4}e^{-2\tilde{\tilde{ x}}} \partial_{\tilde{\tilde{x}}}.\\
\end{array}
\end{equation}
Finally, the transformation
$\hat{x}=\tilde{\tilde{y}},\hat{y}=-\frac{1}{2}
e^{-2\tilde{\tilde{x}}}$ transforms $X,Y$ and $Z$ to the
polynomial form
$$ X=\partial_{\hat{y}}, Y=-{\hat{y}}^2\partial_{\hat{y}}, Z=- 2\hat{y}\partial_{\hat{y}}.$$
%%%%%%%%%%%%%%%%%%%%%%%%%%%%%%%%%%%%%%%%%%%%%%%%%%%%%%%%%%%%%%%%%%%%%%%%%%%%%%%%%%%%%%%%%%%%%%%%%%%%%%
\subsubsection{Realizations of $so(3)$ as vector fields in ${\mathbb{R}}^2$}
Let $L$ be a three-dimensional algebra of vector fields with
negative definite Killing form. Then $L$ has a basis $\{X, Y, Z\}$
such that
$$[X,Y]=Z,[Y,Z]=X,[Z,X]=Y.$$
Find coordinates in which $X=\partial_x$. Then necessarily
$Y-\sqrt{-1}~Z$ has the following form
\begin{equation}\label{te0}
\begin{array}{cc}
Y-\sqrt{-1}~Z=e^{\sqrt{-1} x}\left[ \left( f_1(y)\partial_{x}+f_2(y)\partial_{y}\right) + \sqrt{-1} \left( f_3(y)\partial_{x}+f_4(y)\partial_{y}\right)\right].\\
\end{array}
\end{equation}
Since the rank of $\{\partial_x,
f_1(y)\partial_{x}+f_2(y)\partial_{y},
f_3(y)\partial_{x}+f_4(y)\partial_{y}\}$ cannot be 1, we can
assume without loss of generality that rank of $\{\partial_x,
f_1(y)\partial_{x}+f_2(y)\partial_{y}\}$ is 2. We may therefore
assume that  $f_2(y)\ne 0$.
As $\partial_x, f_1(y)\partial_{x}+f_2(y)\partial_{y}$ are
commuting vectors, there is a change of variables in which
$\partial_x=\partial_{\tilde{x}}$ and
$f_1(y)\partial_{x}+f_2(y)\partial_{y}=\partial_{\tilde{y}}$.
Using section 2, such a change of variables can be given explicitly
as follows:
\begin{equation}
\begin{array}{lll}
\tilde{x}=x-\int{\frac{f_1}{f_2}dy},& \tilde{y}&=\int{\frac{1}{f_2}dy},\\
\end{array}
\end{equation}
Therefore in these coordinates
\begin{equation}\label{te1}
\begin{array}{cc}
Y-\sqrt{-1}~Z=e^{\sqrt{-1}(\tilde{x}+f(\tilde{y}))}\left[\partial_{\tilde{y}} + \sqrt{-1} \left( A(\tilde{y})\partial_{\tilde{x}}+B(\tilde{y})\partial_{\tilde{y}}\right)\right],\\
\end{array}
\end{equation}
where $f(\tilde{y})=\int{\frac{f_1}{f_2}dy},~A(\tilde{y})=
\frac{f_2f_3-f_1f_4}{f_2},$ and $B(\tilde{y})=\frac{f_4}{f_2}.$
Using the fact that $[Y,Z]=\partial_{\tilde{x}}$ if and only if
$[Y-\sqrt{-1}~Z,Y+\sqrt{-1}~Z]=2~\sqrt{-1}\partial_{\tilde{x}}$, the
necessary and sufficient conditions for $f(\tilde{y}), A(\tilde{y})$
and $B(\tilde{y})$ to give a representation of $so(3)$ are
\begin{equation}\label{te2}
\begin{array}{cc}
A^2+ABf'+A'=-1,\\
(AB+B')+(1+B^2)f'=0.\\
\end{array}
\end{equation}
To reduce the form (\ref{te1}) to the simplest form, we look at
the classical Bianchi representation of vector fields on
$\mathbb{P}^2(\mathbb{R})$ induced by the rotations on
$\mathbb{R}^3$. It is given by
$$L_{3:9}:~X=\partial_{\tilde{\tilde{x}}}, Y=\tilde{\tilde{y}}\sin
\tilde{\tilde{x}}
\partial_{\tilde{\tilde{x}}}+(1+{\tilde{\tilde{y}}}^2)\cos \tilde{\tilde{x}} \partial_{\tilde{\tilde{y}}}, Z=\tilde{\tilde{y}}\cos \tilde{\tilde{x}}
\partial_{\tilde{\tilde{x}}}-(1+{\tilde{\tilde{y}}}^2)\sin \tilde{\tilde{x }}\partial_{\tilde{\tilde{y}}}.$$
So
\begin{equation}\label{te3}
\begin{array}{cc}
Y-\sqrt{-1}~Z=e^{\sqrt{-1}\tilde{\tilde{x}}}\left[(1+{\tilde{\tilde{y}}}^2)\partial_{\tilde{\tilde{y}}}- \sqrt{-1}\tilde{\tilde{y}}\partial_{\tilde{\tilde{x}}}\right].\\
\end{array}
\end{equation}
The conditions that the form (\ref{te1}) can be written in the
form (\ref{te3}) with
$\partial{\tilde{\tilde{x}}}=\partial{\tilde{x}}$ and $\tilde{
\tilde{y}}=\psi(\tilde{y})$ are exactly the equations (\ref{te2}).
This gives the transformation
\begin{equation}
\begin{array}{lll}
\tilde{\tilde{x}}=\tilde{x}+f+{\tan}^{-1}B,& \tilde{\tilde{y}}&=-\frac{A}{\sqrt{1+B^2}}.\\
\end{array}
\end{equation}
So the transformation
\begin{equation}
\begin{array}{lll}\label{te4}
\tilde{\tilde{x}}=x+{\tan}^{-1}\left(\frac{f_4}{f_2}\right),&\tilde{ \tilde{y}}&=\frac{f_1f_4-f_2f_3}{\sqrt{f_2^2+f_4^2}}.\\
\end{array}
\end{equation}
maps the form (\ref{te0}) to the form (\ref{te3}).

In case $f_2(y)=0$, $f_4(y)\ne 0$, the formula becomes
\begin{equation}
\begin{array}{lll}\label{te5}
\tilde{\tilde{x}}=x-{\tan}^{-1}\left(\frac{f_2}{f_4}\right)+\frac{\pi}{2},&\tilde{ \tilde{y}}&=\frac{f_1f_4-f_2f_3}{\sqrt{f_2^2+f_4^2}}.\\
\end{array}
\end{equation}
Hence, up to change of coordinates, there is only one realization
as vector fields in ${\mathbb{R}}^2$.
%%%%%%%%%%%%%%%%%%%%%%%%%%%%%%%%%%%%%%%%%%%%%%%%%%%%%%%%%%%%%%%%%%%%%%%%%%%%%%%%%%%%%%%%%%%%%%
\section{Summary of the results}

Based on the discussion in the previous section and using the
notations in ref. \cite{Mahomed1989}, we can state the following
theorem:
\begin{theorem}\label{th1}
Every three-dimensional Lie algebra has one of the following
$17$ realizations in $\mathbb{R}^2$:\\
A) {\rm dim} $G'=0$:\\
Then $G$ is abelian of rank 1 and there are infinitely many
realizations:\\
$L_{3;1}:~X=\partial_y, Y=x\partial_y,Z=f(x)\partial_{y}$ where $f(x)$ is linearly independent of $\{1,x\}$.\\
B) {\rm dim} $(G') =2$:\\
Then the eigenvalues of $G/G'$ on $G'$ never zero and there are
the following eight cases:\\
1) {\rm rank} $(G')=1$ and the eigenvalues of $G/G'$ on $G'$ are real and distinct.\\
$L^{II}_{3:6}:~X=\partial_y, Y=x\partial_y,Z=(c-1)x\partial_{x}-y\partial_{y}, c\neq 0, 1$\\
with the nonzero commutators $[Z,X]=X,[Z,Y]=cY, c\neq 0,1.$\\
2) {\rm rank} $(G')=2$ and the eigenvalues of $G/G'$ on $G'$ are real and distinct.\\
$L^{I}_{3:6}:~X=\partial_x, Y=\partial_y,Z=-x\partial_{x}-cy\partial_{y}, c\neq 0, 1$\\
with the nonzero commutators $[Z,X]=X,[Z,Y]=cY, c\neq 0,1.$\\
3) {\rm rank} $(G')=1$ and the eigenvalues of $G/G'$ on $G'$ are real and repeated with eigenspace of dimension 2.\\
$L^{II}_{3:5}:~X=\partial_y, Y=x\partial_y,Z=-y\partial_{y}$\\
with the nonzero commutators $[Z,X]=X,[Z,Y]=Y$.\\
4) {\rm rank} $(G')=2$ and the eigenvalues of $G/G'$ on $G'$ are real and repeated with eigenspace of dimension 2.\\
$L^{I}_{3:5}:~X=\partial_x, Y=\partial_y,Z=-x\partial_{x}-y\partial_{y}$\\
with the nonzero commutators $[Z,X]=X,[Z,Y]=Y$.\\
5) {\rm rank} $(G')=1$ and the eigenvalues of $G/G'$ on $G'$ are  real and repeated with eigenspace of dimension 1.\\
$L^{II}_{3:3}:~X=x\partial_y, Y=\partial_y,Z=\partial_x-y\,\partial_y$\\
with the nonzero commutators $[Z,X]=X+Y, [Z,Y]=Y$\\
6) {\rm rank} $(G')=2$ and the eigenvalues of $G/G'$ on $G'$ are  real and repeated with eigenspace of dimension 1.\\
$L^{I}_{3:3}:~X=\partial_x, Y=\partial_y,Z=-x\partial_x-(x+y)\,\partial_y$\\
with the nonzero commutators $[Z,X]=X+Y, [Z,Y]=Y$\\
7) {\rm rank} $(G')=1$ and the eigenvalues of $G/G'$ on $G'$ are complex.\\
$L^{II}_{3:7}:~X=\partial_y, Y=x\partial_y, Z=\sin \theta(1+x^2)\, \partial_{x}+y(x\sin \theta- \,\cos \theta)\,\partial_{y}$\\
with the nonzero commutators $[Z,X]=\cos \theta~X-\sin \theta~Y,[Z,Y]=\sin \theta~X+ \cos \theta~Y$.\\
8) {\rm rank} $(G')=2$ and the eigenvalues of $G/G'$ on $G'$ are complex.\\
$L^{I}_{3:7}:~X=\partial_x, Y=\partial_y, Z=(-x \cos \theta -y \sin \theta)\, \partial_{x}+(x\sin \theta-y \,\cos \theta)\,\partial_{y}$\\
with the nonzero commutators $[Z,X]=\cos \theta~X-\sin \theta~Y,[Z,Y]=\sin \theta~X+ \cos \theta~Y$.\\
C)  {\rm dim} $(G')=1$\\
1) The centralizer  $Z_G(G')$ is 2 dimensional, then there are the
two cases:\\
(i) {\rm rank} $(Z_G(G'))$  is 1.\\
$L^{II}_{3:4}:~X=-x\partial_x-y\partial_y, Y=x\partial_y,Z=\partial_y$ \\
with the nonzero commutator $[X,Z]=Z$.\\
(ii) {\rm rank} $(Z_G(G'))$ is 2.\\
$L^{I}_{3:4}:~X=-y\partial_y, Y=\partial_x,Z=\partial_y$ \\
with the nonzero commutator $[X,Z]=Z$.\\
2) The centralizer  $Z_G(G')$ is three-dimensional:\\
$L_{3:2}:~X=\partial_x, Y=x\partial_y,Z=\partial_y$\\
with the nonzero commutator $[X,Y]=Z$.\\
D) $G'=G$:\\
1) The Killing form is negative definite.\\
The Lie algebra is $so(3)$ and there is one realization:\\
$L_{3:9}:~X=\partial_x, Y=y\sin x \partial_x+(1+y^2)\cos x \partial_y, Z=y\cos x \partial_x-(1+y^2)\sin x \partial_y$\\
with the nonzero commutators $[X,Y]=Z,[Y,Z]=X,[Z,X]=Y$.\\
2) The Killing form is indefinite.\\
The Lie algebra is $sl(2,R)$ and there are two cases:\\
(i) The rank of the generators is 2 and there is one of the following three realizations.\\
$L^{I}_{3:8}:~ X=\partial_y, Y=-2xy\partial_x-y^2\partial_y, Z= -2x\partial_x-2y\partial_y,$\\
$L^{II}_{3:8}:~ X=\partial_y, Y=-2xy\partial_x+(-y^2+x^2)\partial_y, Z= -2x\partial_x-2y\partial_y,$\\
$L^{III}_{3:8}:~X=\partial_y, Y=-2xy\partial_x-(y^2+x^2)\partial_y, Z= -2x\partial_x-2y\partial_y,$\\
with the nonzero commutators $[Z,X]=2X,[Z,Y]=-2Y,[X,Y]=Z$.\\
(ii) The rank of the generators is 1 and there is one realization.\\
$L^{IV}_{3:8}:~X=\partial_y, Y=-y^2\partial_y, Z=- 2y\partial_y,$\\
with the nonzero commutators $[Z,X]=2X,[Z,Y]=-2Y,[X,Y]=Z$.\\
\end{theorem}
%%%%%%%%%%%%%%%%%%%%%%%%%%%%%%%%%%%%%%%%%%%%%%%%%%%%%%%%%%%%%%%%%%%%%%%%%%%%%%%%%%%%%%%%%%%%%%
\section{Illustrative examples on the 17 Bianchi types}

The equations considered here have been obtained by determining joint invariants of appropriate 
order for all the Lie-Bianchi types  and transforming the equations by simple point transformations. 
The point of these examples is to recover  something close to the inverse of these 
transformations algorithmically. 

By computing the joint invariants of the realizations of Bianchi
types, we see that there are no second-order invariant ODEs  when
rank $G'=1$. However, there are higher-order invariant ODEs for
such cases. Even when rank $G'=2$, not all of Bianchi types have
second-order invariant ODEs. For this reason, we give a procedure
illustrated by examples given below for each of the Bianchi types
which works in principle for any ODE of arbitrary order that
admits a three-dimensional symmetry algebra to reduce it to its
canonical form.
\begin{example}\rm {$L_{3:1}$} \label{ex1}\\
 Consider the ODE
\begin{equation}\label{1e1}
v^{(4)}=\frac{1}{{v'}^{5}}\left({-v{v'}^{10}+10\,{v'}^{4}v''v'''-15\,{v'}^{3}{v''}^{3}-{v'}^{2}{v'''}^{2}+6\,v'{v''}^{2}v'''-9\,{v''}^{4}}\right)\\
\end{equation}
that admits the three-dimensional point symmetry algebra generated
by
\begin{equation}\label{1e2}
\begin{array}{lll}
Y_1=\frac{\partial}{\partial u},&Y_2=v\frac{\partial}{\partial u},&Y_3=v^2\frac{\partial}{\partial u},\\
\end{array}
\end{equation}
Since $G$ is abelian of rank 1, using Theorem \ref{th1}, the
fourth-order ODE (\ref{1e1}) can be transformed to the canonical
form of $L_{3:1}$ via a point transformation.

In order to construct such a point transformation, one needs to
match the the symmetries with the realizations of the Lie algebra
of $L_{3:1}$ given by Theorem \ref{th1} in the following way:
\begin{equation}\label{1e5}
\begin{array}{ccc}
X=Y_1,&Y=Y_2,&Z=Y_3.\\
\end{array}
\end{equation}
Applying this correspondence to the point transformation
$u=\phi(x,y), v=\psi(x,y)$ yields the system
\begin{equation}\label{1e6}
\begin{array}{cccccc}
\phi_y=1,&x\phi_y=\psi,&f(x)\phi_y={\psi}^2,&\psi_y=0,&x\psi_y=0,&f(x)\psi_y=0.\\
\end{array}
\end{equation}
The solution of the system (\ref{1e6}) gives the following point
transformation
\begin{equation}\label{1e7}
\begin{array}{cc}
u=y,& v=x,\\
\end{array}
\end{equation}
for $f(x)=x^2$ which transforms ODE (\ref{1e1}) to its canonical
form
\begin{equation}\label{1e8}
y^{(4)}=\left(\frac{f^{(4)}}{f''}\right)y''+g\left(x,y'''-\left(\frac{f'''}{f''}\right)y''\right),\\
\end{equation}
with $g(z,w)=z+w^2$.
\end{example}
%%%%%%%%%%%%%%%%%%%%%%%%%%%%%%%%%%%%%%%%%%%%%%%%%%%%%%%%%%%%%%%%%%%%%%%%%%%%%%%%%%%%%%%%%%%%%%%%%%%%%%%%%%%%%%%%%%%%%%%%%%%%%%
\begin{example}\rm {$L_{3:2}$} \label{ex2}\\
 Consider the ODE
\begin{equation}\label{2e1}
v'''=-\left( v'-v'' \right) ^{3}{{\rm e}^{- 3\,u}}+{{\rm e}^{3\,u}}-  2\,v'+3\,v'',\\
\end{equation}
that admits the three-dimensional point symmetry
algebra generated by
\begin{equation}\label{2e2}
\begin{array}{lll}
Y_1=\frac{\partial}{\partial v},&Y_2=e^{u}\frac{\partial}{\partial v},&Y_3=e^{-u}\frac{\partial}{\partial u},\\
\end{array}
\end{equation}
with the nonzero commutators
\begin{equation}\label{2e3}
\begin{tabular}{llllll}
$[Y_{2}, Y_{3}] =-Y_{1}$.\\
\end{tabular}
\end{equation}
Since dim $G'=1$, dim $Z_G(G')=3$, using Theorem \ref{th1}, the
third-order ODE (\ref{2e1}) can be transformed to the canonical
form of $L_{3:2}$ via a point transformation.

In order to construct such a point transformation, one needs to
match  any vector from $G'=<Y_1>$ with $Z=\frac{\partial
}{\partial y}$ and any vector  which is functionally independent
of $G'$ with $X=\frac{\partial }{\partial x}$ in the following
way:
\begin{equation}\label{2e5}
\begin{array}{ccc}
Z=Y_1,&X=Y_3.\\
\end{array}
\end{equation}
Applying this correspondence to the point transformation
$u=\phi(x,y), v=\psi(x,y)$ yields the system
\begin{equation}\label{2e6}
\begin{array}{cccc}
\phi_x=e^{-\phi},&\phi_y=0,&\psi_x=0,&\psi_y=1.\\
\end{array}
\end{equation}
The solution of the system (\ref{2e6}) gives the following point
transformation
\begin{equation}\label{2e7}
\begin{array}{cc}
u=\ln{x},& v=y,\\
\end{array}
\end{equation}
which transforms ODE (\ref{2e1}) to its canonical form
\begin{equation}\label{2e8}
y'''=f\left(y''\right),
\end{equation}
with $f(z)=z^3+1$.
\end{example}
%%%%%%%%%%%%%%%%%%%%%%%%%%%%%%%%%%%%%%%%%%%%%%%%%%%%%%%%%%%%%%%%%%%%%%%%%%%%%%%%%%%%%%%%%%%%%%%%%%%%%%%%%%%%%%%%%%%%%%%%%%%%%%
\begin{example}\rm {$L^{I}_{3:3}$} \label{ex3}\\
 Consider the ODE
\begin{equation}\label{3e1}
v''=-\frac{1}{3}(v'-2)^3\exp{\left(\frac{v'+1}{v'-2}\right)}\\
\end{equation}
that admits the three-dimensional point symmetry algebra generated
by
\begin{equation}\label{3e2}
\begin{array}{lll}
Y_1=\frac{\partial}{\partial v},&Y_2=\frac{\partial}{\partial u},&Y_3=(5u-v)\frac{\partial}{\partial u}+(4u+v)\frac{\partial}{\partial v},\\
\end{array}
\end{equation}
with the nonzero commutators
\begin{equation}\label{3e3}
\begin{tabular}{llllll}
$[Y_{1}, Y_{3}] = Y_{1}-Y_{2}$,& $[Y_{2}, Y_{3}] = 4Y_{1}+5Y_{2}$.\\
\end{tabular}
\end{equation}
Here dim $G'=2$, rank $(G')=2$ and the adjoint action of
$G/G'=<\overline{Y}_3>$ on $G'=<Y_1,Y_2>$ given by
\begin{equation}\label{3e4}
{\rm ad}\, (\overline{Y}_3)= \left({
\begin{array}{cc}
-1 & -4 \\
1  & -5 \\
\end{array}
}\right)
\end{equation}
has $\lambda=-3$ as a repeated real eigenvalue. The vector
$2Y_1+Y_2$ is an eigenvector and $Y_1$ is a generalized
eigenvector because, in two dimensions, any vector linearly
independent of the eigenvector is a generalized eigenvector. Using
Theorem \ref{th1}, the second-order ODE (\ref{3e1}) can be
transformed to the canonical form of $L^{I}_{3:3}$ via a point
transformation.

In order to construct such a point transformation, one needs to
match the the generalized eigenvector with $X=\frac{\partial
}{\partial x}$ and the scaled eigenvector by $\frac{1}{\lambda}$
with $Y=\frac{\partial }{\partial y}$, in the following way:
\begin{equation}\label{3e5}
\begin{array}{ccc}
X=Y_1,&Y=-\frac{1}{3}(2Y_1+Y_2).\\
\end{array}
\end{equation}
Applying this correspondence to the point transformation
$u=\phi(x,y), v=\psi(x,y)$ yields the system
\begin{equation}\label{3e6}
\begin{array}{cccc}
\phi_x=0,&\phi_y=-\frac{1}{3},&\psi_x=1,&\psi_y=-\frac{2}{3}.\\
\end{array}
\end{equation}
The solution of the system (\ref{3e6}) gives the following point
transformation
\begin{equation}\label{3e7}
\begin{array}{cc}
u=-\frac{1}{3}y,& v=x-\frac{2}{3}y,\\
\end{array}
\end{equation}
which transforms ODE (\ref{3e1}) to its canonical form
\begin{equation}\label{3e8}
y''=C \exp{(-y')}
\end{equation}
with $C=-e$.
\end{example}
%%%%%%%%%%%%%%%%%%%%%%%%%%%%%%%%%%%%%%%%%%%%%%%%%%%%%%%%%%%%%%%%%%%%%%%%%%%%%%%%%%%%%%%%%%%%%%%%%%%%%%%%%%%%%%%%%%%%%%%%%%%%%%
\begin{example}\rm {$L^{II}_{3:3}$} \label{ex4}\\
 Consider the ODE
\begin{equation}\label{4e1}
v'''=\frac{1}{u}\left({{\rm e}^{-u}}\ln  \left(  \left( uv''+2\,v'\right) { {\rm e}^{u}}-u \right) -{{\rm e}^{-u}}u-3\,v''\right)\\
\end{equation}
that admits the three-dimensional point symmetry algebra generated
by
\begin{equation}\label{4e2}
\begin{array}{lll}
Y_1=\frac{\partial}{\partial v},&Y_2=\frac{1}{u}\frac{\partial}{\partial v},&Y_3=\frac{\partial}{\partial u}+(\frac{e^{-u}-v-uv}{u})\frac{\partial}{\partial v},\\
\end{array}
\end{equation}
with the nonzero commutators
\begin{equation}\label{4e3}
\begin{tabular}{llllll}
$[Y_{1}, Y_{3}] = -Y_{2}-Y_{1}$,& $[Y_{2}, Y_{3}] = -Y_{2}$.\\
\end{tabular}
\end{equation}
Here dim $G'=2$, rank $(G')=1$ and the adjoint action of
$G/G'=<\overline{Y}_3>$ on $G'=<Y_1,Y_2>$ given by
\begin{equation}\label{4e4}
{\rm ad}\,(\overline{Y}_3)= \left({
\begin{array}{cc}
1 & 0 \\
1  & 1 \\
\end{array}
}\right)
\end{equation}
has $\lambda=1$ as a repeated real eigenvalue. Here $Y_2$ is an
eigenvector and $Y_1$ is a generalized eigenvector  for the same
reason as explained in example 6.3. Using the Theorem \ref{th1},
the second-order ODE (\ref{4e1}) can be transformed to the
canonical form of $L^{II}_{3:3}$ via a point transformation.

In order to construct such a point transformation, one needs to
match the the generalized eigenvector with $X=x\frac{\partial
}{\partial y}$ and the eigenvector scaled  by $\frac{1}{\lambda}$
with $Y=\frac{\partial }{\partial y}$ in the following way:
\begin{equation}\label{4e5}
\begin{array}{ccc}
X=Y_1,&Y=Y_2.\\
\end{array}
\end{equation}
Applying this correspondence to the point transformation
$u=\phi(x,y), v=\psi(x,y)$ yields the system
\begin{equation}\label{4e6}
\begin{array}{cccc}
\phi_y=0,&x\phi_y=0,&\psi_y=\frac{1}{\phi},&x\psi_y=1.\\
\end{array}
\end{equation}
The solution of the system (\ref{4e6}) gives the following point
transformation
\begin{equation}\label{4e7}
\begin{array}{cc}
u=x,& v=\frac{y}{x},\\
\end{array}
\end{equation}
which transforms the vector $Y_3$ which is linearly independent of
$G'$ to
\begin{equation}\label{4e8}
\begin{array}{cc}
Z=\frac{\partial}{\partial x}+(-y+f(x))\frac{\partial}{\partial y},\\
\end{array}
\end{equation}
with $f(x)=e^{-x}$. Such a function $f(x)$ can be absorbed using
the transformation
\begin{equation}\label{4e9}
\begin{array}{cc}
\tilde{x}=x,& \tilde{y}=y-e^{-x}\int e^x f(x) dx=y-xe^{-x}.\\
\end{array}
\end{equation}
Finally, the composition of the transformations (\ref{4e7}) and
(\ref{4e9}) transforms ODE (\ref{4e1}) to its canonical form
\begin{equation}\label{4e10}
\tilde{y}'''=e^{-\tilde{x}}g\left(e^{\tilde{x}}\tilde{y}''\right),
\end{equation}
with $g(z)=-3+\ln (z-2)$.
\end{example}
%%%%%%%%%%%%%%%%%%%%%%%%%%%%%%%%%%%%%%%%%%%%%%%%%%%%%%%%%%%%%%%%%%%%%%%%%%%%%%%%%%%%%%%%%%%%%%%%%%%%%%%%%%%%%%%%%%%%%%%%%%%%%%
\begin{example}\rm {$L^{I}_{3:4}$} \label{ex5}\\
 Consider the ODE
\begin{equation}\label{5e1}
v'''=-{\frac {  {v'}^{4}-3\,{v''}^{2}  { }}{{v'} }}-{\frac { { {\rm e}^{-\,v}} \left({v'}^{2}-v'' \right) ^{3}  { }}{{v'}^{2} \left( uv'+1 \right) ^{3}}}\\
\end{equation}
that admits the three-dimensional point symmetry algebra generated
by
\begin{equation}\label{5e2}
\begin{array}{lll}
Y_1=-u\frac{\partial}{\partial u}+\frac{\partial}{\partial v},&Y_2=e^{-v}\frac{\partial}{\partial u},&Y_3=-ue^{-v}\frac{\partial}{\partial u}+e^{-v}\frac{\partial}{\partial v},\\
\end{array}
\end{equation}
with the nonzero commutators
\begin{equation}\label{5e3}
\begin{tabular}{llllll}
$[Y_{1}, Y_{3}] =-Y_{3}$.\\
\end{tabular}
\end{equation}
Since dim $G'=1$, dim $Z_G(G')=2$ and rank $(Z_G(G'))=2$, using
Theorem \ref{th1}, the third-order ODE (\ref{5e1}) can be
transformed to the canonical form of $L^{I}_{3:4}$ via a point
transformation.

In order to construct such a point transformation, one needs to
match  any vector from $G'=<Y_3>$ with $Z=\frac{\partial
}{\partial y}$ and any vector from $Z_G(G')=<Y_3,Y_2>$ which is
linearly independent of $G'$ with $Y=\frac{\partial }{\partial x}$
in the following way:
\begin{equation}\label{5e5}
\begin{array}{ccc}
Z=Y_3,&Y=Y_2.\\
\end{array}
\end{equation}
Applying this correspondence to the point transformation
$u=\phi(x,y), v=\psi(x,y)$ yields the system
\begin{equation}\label{5e6}
\begin{array}{cccc}
\phi_x=e^{-\psi},&\phi_y=-\phi e^{-\psi},&\psi_x=0,&\psi_y=e^{-\psi}.\\
\end{array}
\end{equation}
The solution of the system (\ref{5e6}) gives the following point
transformation
\begin{equation}\label{5e7}
\begin{array}{cc}
u=\frac{x}{y},& v=\ln{y},\\
\end{array}
\end{equation}
which transforms ODE (\ref{5e1}) to its canonical form
\begin{equation}\label{5e8}
y'''={y'}f\left(\frac{y''}{{y'}}\right),
\end{equation}
with $f(z)=z^3+3z^2$.
\end{example}
%%%%%%%%%%%%%%%%%%%%%%%%%%%%%%%%%%%%%%%%%%%%%%%%%%%%%%%%%%%%%%%%%%%%%%%%%%%%%%%%%%%%%%%%%%%%%%%%%%%%%%%%%%%%%%%%%%%%%%%%%%%%%%
\begin{example}\rm {$L^{II}_{3:4}$} \label{ex6}\\
Consider the ODE
\begin{equation}\label{6e1}
v'''={\frac { \left( 32\,{v}^{6}-1 \right){v'}^{8}-8\,{v}^{4}{v'}^{5}v''+12\,{v}^{5}{v'}^{3}{v''}^{2}-3\,{v}^{3}{v''}^{3}}{4\,{v}^{5}{v'}^{4}-{v}^{3}v'v''}}\\
\end{equation}
that admits the three-dimensional point symmetry algebra generated
by
\begin{equation}\label{6e2}
\begin{array}{lll}
Y_1=\frac{\partial}{\partial u},&Y_2=v\frac{\partial}{\partial u},&Y_3=(u-v^4)\frac{\partial}{\partial u}+v\frac{\partial}{\partial v},\\
\end{array}
\end{equation}
with the nonzero commutators
\begin{equation}\label{6e3}
\begin{tabular}{llllll}
$[Y_{1}, Y_{3}] =Y_{1}$.\\
\end{tabular}
\end{equation}
Since dim $G'=1$, dim $Z_G(G')=2$ and rank $(Z_G(G'))=1$, using
Theorem \ref{th1}, the third-order ODE (\ref{6e1}) can be
transformed to the canonical form of $L^{II}_{3:4}$ via a point
transformation.

In order to construct such a point transformation, one needs to
match  any vector from $G'=<Y_1>$ with $Z=\frac{\partial
}{\partial y}$ and any vector from $Z_G(G')=<Y_1,Y_2>$ which is
linearly independent of $G'$ with $Y=x\frac{\partial }{\partial
y}$, in the following way:
\begin{equation}\label{6e5}
\begin{array}{ccc}
Z=Y_1,&Y=Y_2.\\
\end{array}
\end{equation}
Applying this correspondence to the point transformation
$u=\phi(x,y), v=\psi(x,y)$ yields the system
\begin{equation}\label{6e6}
\begin{array}{cccc}
\phi_y=1,&x\phi_y=\psi,&\psi_y=0,&x\psi_y=0.\\
\end{array}
\end{equation}
The solution of the system (\ref{6e6}) gives the following point
transformation
\begin{equation}\label{6e7}
\begin{array}{cc}
u=y,& v=x,\\
\end{array}
\end{equation}
which transforms the vector $-Y_3$ which is linearly independent
of $Z_G(G')$ to
\begin{equation}\label{6e8}
\begin{array}{cc}
X=-x\frac{\partial}{\partial x}+(-y+f(x))\frac{\partial}{\partial y},\\
\end{array}
\end{equation}
with $f(x)=x^4$. Such a function $f(x)$ can be absorbed using the
transformation
\begin{equation}\label{6e9}
\begin{array}{cc}
\tilde{x}=x,& \tilde{y}=y+x\int \frac{f(x)}{x^2}dx=y+\frac{1}{3}x^4.\\
\end{array}
\end{equation}
Finally, the composition of the transformations (\ref{6e7}) and
(\ref{6e9}) transforms ODE (\ref{6e1}) to its canonical form
\begin{equation}\label{6e10}
\tilde{y}'''=\frac{1}{\tilde{x}^2}g\left(\tilde{x}\tilde{y}''\right),
\end{equation}
with $g(z)=\frac{1}{z}$.
\end{example}
%%%%%%%%%%%%%%%%%%%%%%%%%%%%%%%%%%%%%%%%%%%%%%%%%%%%%%%%%%%%%%%%%%%%%%%%%%%%%%%%%%%%%%%%%%%%%%%%%%%%%%%%%%%%%%%%%%%%%%%%%%%%%%
\begin{example}\rm {$L^{I}_{3:5}$} \label{ex7}\\
Consider the ODE
\begin{equation}\label{7e1}
v'''= \frac{1}{{{v'}^{4}}}\left({- \left( {v''}^{2}+v'' \right) ^{2}{{\rm e}^{-3\,v}}+2\,{v'}^{7}+3\,v''{v'}^{5}+3\,{v'}^{3}{v''}^{2}}\right)\\
\end{equation}
that admits the three-dimensional point symmetry algebra generated
by
\begin{equation}\label{7e2}
\begin{array}{lll}
Y_1=\frac{\partial}{\partial u},&Y_2=u\frac{\partial}{\partial u}+\frac{\partial}{\partial v},&Y_3=e^{-v}\frac{\partial}{\partial v},\\
\end{array}
\end{equation}
with the nonzero commutators
\begin{equation}\label{7e3}
\begin{tabular}{llllll}
$[Y_{1}, Y_{2}] = Y_{1}$,& $[Y_{2}, Y_{3}] = -Y_{3}$.\\
\end{tabular}
\end{equation}
Here dim $G'=2$, rank $(G')=2$ and the adjoint action of
$G/G'=<\overline{Y}_2>$ on $G'=<Y_1,Y_3>$  is given by
\begin{equation}\label{7e4}
{\rm ad}\,(\overline{Y}_2)= \left({
\begin{array}{cc}
-1 & 0 \\
0  & -1 \\
\end{array}
}\right).
\end{equation}
So $\lambda=-1$ is a repeated real eigenvalue with eigenspace of
dimension 2. Using Theorem \ref{th1}, the third-order ODE
(\ref{7e1}) can be transformed to the canonical form of
$L^{I}_{3:5}$ via a point transformation.

In order to construct such a point transformation, one needs to
match any two linearly independent vectors of $G'$ with
$X=\frac{\partial }{\partial x}$ and $Y=\frac{\partial }{\partial
y}$. For example, one can try the obvious choices:
\begin{equation}\label{7e5}
\begin{array}{ccc}
X=Y_1,&Y=Y_3,\\
\end{array}
\end{equation}
or the opposite
\begin{equation}\label{7e6}
\begin{array}{ccc}
X=Y_3,&Y=Y_1.\\
\end{array}
\end{equation}
Applying the correspondence (\ref{7e5}) to the point
transformation $u=\phi(x,y), v=\psi(x,y)$ yields the system
\begin{equation}\label{7e7}
\begin{array}{cccc}
\phi_x=1,&\phi_y=0,&\psi_x=0,&\psi_y=e^{-\psi}.\\
\end{array}
\end{equation}
The solution of the system (\ref{7e7}) gives the following point
transformation
\begin{equation}\label{7e8}
\begin{array}{cc}
u=x,& v=\ln{y},\\
\end{array}
\end{equation}
which transforms ODE (\ref{7e1}) to its canonical form
\begin{equation}\label{7e9}
y'''=f(y') {y''}^2
\end{equation}
with $f(z)=\frac{3z^3-1}{z^4}$.

Similarly, applying the correspondence (\ref{7e6}) to the point
transformation $u=\phi(x,y), v=\psi(x,y)$ yields the system
\begin{equation}\label{7e10}
\begin{array}{cccc}
\phi_x=0,&\phi_y=1,&\psi_x=e^{-\psi},&\psi_y=0.\\
\end{array}
\end{equation}
The solution of the system (\ref{7e10}) gives the following point
transformation
\begin{equation}\label{7e11}
\begin{array}{cc}
u=y,& v=\ln{x},\\
\end{array}
\end{equation}
which transforms ODE (\ref{7e1}) to its canonical form
\begin{equation}\label{7e12}
y'''=f(y') {y''}^2
\end{equation}
with $f(z)=z^2$.
\end{example}
%%%%%%%%%%%%%%%%%%%%%%%%%%%%%%%%%%%%%%%%%%%%%%%%%%%%%%%%%%%%%%%%%%%%%%%%%%%%%%%%%%%%%%%%%%%%%%%%%%%%%%%%%%%%%%%%%%%%%%%%%%%%%%
\begin{example}\rm {$L^{II}_{3:5}$} \label{ex8}\\
Consider the ODE
\begin{equation}\label{8e1}
v^{(4)}={\frac {4\,v{v'}^{12}-4\,v{v'}^{9}v''+v{v'}^{6}{v''}^{2}+10\,{v'}^{2}v''{v'''}^{2}-45\,v'{v''}^{3}v'''+45\,{v''}^{5}}{{v'}^{2} \left( v'''v'-3\,{v''}^{2} \right) }}\\
\end{equation} that admits the three-dimensional point
symmetry algebra generated by
\begin{equation}\label{8e2}
\begin{array}{lll}
Y_1=\frac{\partial}{\partial u},&Y_2=v\frac{\partial}{\partial u},&Y_3=(u+v^2)\frac{\partial}{\partial u},\\
\end{array}
\end{equation}
with the nonzero commutators
\begin{equation}\label{8e3}
\begin{tabular}{llllll}
$[Y_{1}, Y_{3}] =Y_{1}$,& $[Y_{2}, Y_{3}] =Y_{2}$.\\
\end{tabular}
\end{equation}
Here dim $G'=2$, rank $(G')=1$ and the adjoint action of
$G/G'=<\overline{Y}_3>$ on \\$G'=<Y_1,Y_2>$ is given by
\begin{equation}\label{8e4}
{\rm ad}\, (\overline{Y}_3)= \left({
\begin{array}{cc}
-1 & 0 \\
0  & -1 \\
\end{array}
}\right).
\end{equation}
 We have $\lambda=-1$ as a repeated real eigenvalue with eigenspace of
dimension 2. Using Theorem \ref{th1}, the forth-order ODE
(\ref{8e1}) can be transformed to the canonical form of
$L^{II}_{3:5}$ via a point transformation.

In order to construct such a point transformation, one needs to
match any two linearly independent vectors of $G'$ with
$X=\frac{\partial }{\partial y}$ and $Y=x\frac{\partial }{\partial
y}$. For example, one can try the obvious choice:
\begin{equation}\label{8e5}
\begin{array}{ccc}
X=Y_1,&Y=Y_2,\\
\end{array}
\end{equation}
Applying the correspondence (\ref{8e5}) to the point
transformation $u=\phi(x,y), v=\psi(x,y)$ yields the system
\begin{equation}\label{8e7}
\begin{array}{cccc}
\phi_y=1,&x\phi_y=\psi,&\psi_y=0,&x\psi_y=0.\\
\end{array}
\end{equation}
The solution of the system (\ref{8e7}) gives the following point
transformation
\begin{equation}\label{8e8}
\begin{array}{cc}
u=y,& v=x,\\
\end{array}
\end{equation}
which transforms the vector $-Y_3$ which is linearly independent
of $G'$ to
\begin{equation}\label{8e9}
\begin{array}{cc}
Z=(-y+f(x))\frac{\partial}{\partial y},\\
\end{array}
\end{equation}
with $f(x)=-x^2$. The function $f(x)$ can be absorbed using the
transformation
\begin{equation}\label{8e10}
\begin{array}{cc}
\tilde{x}=x,& \tilde{y}=y-f(x)=y+x^2.\\
\end{array}
\end{equation}
Finally, the composition of the transformations (\ref{8e8}) and
(\ref{8e10}) transforms ODE (\ref{8e1}) to its canonical form
\begin{equation}\label{8e11}
\tilde{y}^{(4)}=\tilde{y}''~g\left(x,\frac{\tilde{y}'''}{\tilde{y}''}\right),\\
\end{equation}
with $g(z,w)=\frac{z}{w}$.
\end{example}
%%%%%%%%%%%%%%%%%%%%%%%%%%%%%%%%%%%%%%%%%%%%%%%%%%%%%%%%%%%%%%%%%%%%%%%%%%%%%%%%%%%%%%%%%%%%%%%%%%%%%%%%%%%%%%%%%%%%%%%%%%%%%%
\begin{example}\rm {$L^{I}_{3:6}$}\label{ex9}\\
 Consider the ODE
\begin{equation}\label{9e1}
v''=\frac{1}{9}v^{-4}\left(18\,{v_{{}}}^{3}{v'}^{2}-\sqrt{{v_{{}}}^{2}-2\,v'} \left( {v_{{}}}^{2}+v' \right) ^{5/2}\right)\\
\end{equation}
that admits the three-dimensional point symmetry algebra generated
by
\begin{equation}\label{9e2}
\begin{array}{lll}
Y_1=\frac{\partial}{\partial u},&Y_2=v^2\frac{\partial}{\partial v},&Y_3=(u+\frac{4}{5v})\frac{\partial}{\partial u}-(\frac{2}{5}uv^2+\frac{7}{5}v)\frac{\partial}{\partial v},\\
\end{array}
\end{equation}
with the nonzero commutators
\begin{equation}\label{9e3}
\begin{tabular}{llllll}
$[Y_{1}, Y_{3}] = Y_{1}-\frac{2}{5}Y_{2}$,& $[Y_{2}, Y_{3}] = -\frac{4}{5}Y_{1}+\frac{7}{5}Y_{2}$.\\
\end{tabular}
\end{equation}
We have dim $G'=2$, rank $(G')=2$ and the adjoint action of
$G/G'=<\overline{Y}_3>$ on $G'=<Y_1,Y_2>$ given by
\begin{equation}\label{9e4}
{\rm ad}\, (\overline{Y}_3)= \left({
\begin{array}{cc}
-1 & \frac{4}{5} \\
\frac{2}{5}  & -\frac{7}{5} \\
\end{array}
}\right)
\end{equation}
has $\lambda_1=-\frac{3}{5}$ and $\lambda_2=-\frac{9}{5}$ as
distinct real eigenvalues. The corresponding eigenvectors are
$2Y_1+Y_2$ and $Y_2-Y_1$ respectively. Using the Theorem
\ref{th1}, the second-order ODE (\ref{9e1}) can be transformed to
the canonical form of $L^{I}_{3:6}$ via a point transformation.

In order to construct such a point transformation, one needs to
match the the two eigenvectors of ad $(G/G')$ on $G'$ with constant
multiples $X=\frac{\partial }{\partial x}$ and $Y=\frac{\partial
}{\partial y}$ respectively, in the following way:
\begin{equation}\label{9e5}
\begin{array}{ccc}
X=r(2Y_1+Y_2),&Y=s(Y_2-Y_1),& r, s\in \mathbb{R}\setminus{\{0\}}.\\
\end{array}
\end{equation}
Applying the correspondence (\ref{9e5}) to the point
transformation $u=\phi(x,y), v=\psi(x,y)$ yields the system
\begin{equation}\label{9e6}
\begin{array}{cccc}
\phi_x=2r,&\phi_y=-s,&\psi_x=r{\psi}^2,&\psi_y=s{\psi}^2.\\
\end{array}
\end{equation}
The solution of the system (\ref{9e6}) gives the following point
transformation
\begin{equation}\label{9e7}
\begin{array}{cc}
u=2rx-sy,& v=\frac{-1}{rx+sy},\\
\end{array}
\end{equation}
which transforms ODE (\ref{9e1}) to its canonical form
\begin{equation}\label{9e8}
y''=C y'^{\frac{c-2}{c-1}}, c\ne 0,\frac{1}{2},1,2,
\end{equation}
with $C=r^{\frac{3}{2}}{(-s)}^{-\frac{1}{2}}$ and
$c=\frac{\lambda_2}{\lambda_1}=3$.
\end{example}
%%%%%%%%%%%%%%%%%%%%%%%%%%%%%%%%%%%%%%%%%%%%%%%%%%%%%%%%%%%%%%%%%%%%%%%%%%%%%%%%%%%%%%%%%%%%%%%%%%%%%%%%%%%%%%%%%%%%%%%%%%%%%%
\begin{example}\rm {$L^{II}_{3:6}$} \label{ex10}\\
 Consider the ODE
\begin{equation}\label{10e1}
v'''= \frac{1}{2} \left(-4\,v'+6\,v''-8 \right)+ \frac{1}{2}\left(v'-v''+2 \right) ^{ 2} {{\rm e}^{2\,u}}\\
\end{equation}
that admits the three-dimensional point symmetry algebra generated
by
\begin{equation}\label{10e2}
\begin{array}{lll}
Y_1=\frac{\partial}{\partial v},&Y_2=e^u\frac{\partial}{\partial v},&Y_3=\frac{\partial}{\partial u}-(4u+2v)\frac{\partial}{\partial v},\\
\end{array}
\end{equation}
with the nonzero commutators
\begin{equation}\label{10e3}
\begin{tabular}{llllll}
$[Y_{1}, Y_{3}] =-2 Y_{1}$,& $[Y_{2}, Y_{3}] = -3Y_{2}$.\\
\end{tabular}
\end{equation}
Here dim $G'=2$, rank $(G')=1$ and the adjoint action of
$G/G'=<\overline{Y}_3>$ on $G'=<Y_1,Y_2>$  has $\lambda_1=2$ and
$\lambda_2=3$ as distinct real eigenvalues with eigenvector $Y_1$
and $Y_2$ respectively. Using Theorem \ref{th1}, the third-order
ODE (\ref{10e1}) can be transformed to the canonical form of
$L^{II}_{3:6}$  via a point transformation.

In order to construct such a point transformation, one needs to
match the the two eigenvectors of ad $(G/G')$ on $G'$ with
$X=\frac{\partial }{\partial y}$ and $Y=x\frac{\partial }{\partial
y}$ respectively, in the following way:
\begin{equation}\label{10e5}
\begin{array}{ccc}
X=Y_1,&Y=Y_2.\\
\end{array}
\end{equation}
Applying the correspondence (\ref{10e5}) to the point
transformation $u=\phi(x,y), v=\psi(x,y)$ yields the system
\begin{equation}\label{10e6}
\begin{array}{cccc}
\phi_y=0,&x\phi_y=0,&\psi_y=1,&x\psi_y=e^{\phi}.\\
\end{array}
\end{equation}
The solution of the system (\ref{10e6}) gives the following point
transformation
\begin{equation}\label{10e7}
\begin{array}{cc}
u=\ln x,& v=y,\\
\end{array}
\end{equation}
which transforms the vector $\frac{1}{2}Y_3$ which is linearly
independent of $G'$ to
\begin{equation}\label{10e8}
\begin{array}{cc}
Z=(c-1)x\frac{\partial}{\partial x}+(-y+f(x))\frac{\partial}{\partial y},\\
\end{array}
\end{equation}
with $f(x)=-2\ln x$ and
$c=\frac{\lambda_2}{\lambda_1}=\frac{3}{2}$. The function $f(x)$
can be absorbed using the transformation
\begin{equation}\label{10e9}
\begin{array}{cc}
\tilde{x}=x,& \tilde{y}=y+\frac{1}{1-c}~x^{\frac{1}{1-c}}\int f(x) x^{\frac{2-c}{c-1}}dx=y+2\ln x-1.\\
\end{array}
\end{equation}
Finally, the composition of the transformations (\ref{10e7}) and
(\ref{10e9}) transforms ODE (\ref{10e1}) to its canonical form
\begin{equation}\label{10e10}
\tilde{y}'''=\tilde{x}^{\frac{2-3c}{c-1}}g\left(\tilde{y}''\tilde{x}^{\frac{2c-1}{c-1}}\right), c\ne 0,1\\
\end{equation}
with $g(z)=\frac{1}{2}z^2$ and $c=\frac{3}{2}$.
\end{example}
%%%%%%%%%%%%%%%%%%%%%%%%%%%%%%%%%%%%%%%%%%%%%%%%%%%%%%%%%%%%%%%%%%%%%%%%%%%%%%%%%%%%%%%%%%%%%%%%%%%%%%%%%%%%%%%%%%%%%%%%%%%%%%
\begin{example}\rm {$L^{I}_{3:7}$}\label{ex11}\\
 Consider the ODE
\begin{equation}\label{11e1}
v''=\frac{1}{9} {(2{v'}^2-2v'+5)}^{\frac{3}{2}}\exp{\left(3\arctan{\left(\frac{ v'-2}{v'+1}\right)}\right)}\\
\end{equation}
that admits the three-dimensional point symmetry algebra generated
by
\begin{equation}\label{11e2}
\begin{array}{lll}
Y_1=\frac{\partial}{\partial u},&Y_2=\frac{\partial}{\partial v},&Y_3=(4u+v)\frac{\partial}{\partial u}+(5v-\frac{5}{2}u)\frac{\partial}{\partial v},\\
\end{array}
\end{equation}
with the nonzero commutators
\begin{equation}\label{11e3}
\begin{tabular}{llllll}
$[Y_{1}, Y_{3}] =4 Y_{1}-\frac{5}{2}Y_{2}$,& $[Y_{2}, Y_{3}] = Y_{1}+5Y_{2}$.\\
\end{tabular}
\end{equation}
Here dim $G'=2$, rank $(G')=2$ and the adjoint action of
$G/G'=<\overline{Y}_3>$ on $G'=<Y_1,Y_2>$ is given by
\begin{equation}\label{11e4}
{\rm ad}\, (\overline{Y}_3)= \left({
\begin{array}{cc}
-4 & -1 \\
\frac{5}{2}  &-5\\
\end{array}
}\right).
\end{equation}
The eigenvalues are $-\frac{9}{2}\pm\frac{3}{2}\sqrt{-1}$  with
eigenvectors $(\frac{1}{5}\pm\frac{3}{5}\sqrt{-1})Y_1+Y_2$
respectively. Using the Theorem \ref{th1}, the second-order ODE
(\ref{11e1}) can be transformed to the canonical form of
$L^{I}_{3:7}$ via a point transformation.

In order to construct such a point transformation, one needs to
match the real and imaginary parts of an eigenvector of ad $(G/G')$
on $G'$ with $X=\frac{\partial }{\partial x}$ and
$Y=\frac{\partial }{\partial y}$ respectively in the following
way:
\begin{equation}\label{11e5}
\begin{array}{cc}
X=\frac{1}{5}Y_1+Y_2,&Y=\frac{3}{5}Y_1.\\
\end{array}
\end{equation}
Applying the correspondence (\ref{11e5}) to the point
transformation $u=\phi(x,y), v=\psi(x,y)$ yields the system
\begin{equation}\label{11e6}
\begin{array}{cccc}
\phi_x=\frac{1}{5},&\phi_y=\frac{3}{5},&\psi_x=1,&\psi_y=0.\\
\end{array}
\end{equation}
Solving the system (\ref{11e6}) gives the following point
transformation
\begin{equation}\label{11e7}
\begin{array}{cc}
u=\frac{1}{5}x+\frac{3}{5}y,& v=x,\\
\end{array}
\end{equation}
which transforms ODE (\ref{11e1}) after simplification using the
identity $$\tan^{-1} x-\tan^{-1}
\left(\frac{2x-1}{x+2}\right)=c=\tan^{-1}
\left(\frac{1}{2}\right)$$ to its canonical form
\begin{equation}\label{11e8}
y''=C {(1+{y'}^2)}^{\frac{3}{2}}\exp{(b\arctan{ y'})}\\
\end{equation}
 with $C=-\frac{e^{3c}}{\sqrt{5}}$ and $b=\cot{\theta}=-3$ where
 $\theta=arg(-\frac{9}{2}+\frac{3}{2}\sqrt{-1})$.
\end{example}
%%%%%%%%%%%%%%%%%%%%%%%%%%%%%%%%%%%%%%%%%%%%%%%%%%%%%%%%%%%%%%%%%%%%%%%%%%%%%%%%%%%%%%%%%%%%%%%%%%%%%%%%%%%%%%%%%%%%%%%%%%%%%%
\begin{example}\rm {$L^{II}_{3:7}$} \label{ex12}\\
 Consider the ODE
\begin{equation}\label{12e1}
v'''={\frac {4{{\rm e}^{-8\,\arctan \left( \frac{1}{u} \right)}}}{ \left( {u}^{2}+1 \right) ^{\frac{5}{2}}u \left( 1-4 \left( {u}^{2}+1 \right) ^{\frac{3}{2}} \left( uv''+2\,v' \right)  \right)  }}-3\,{\frac {2{u}^{2}v''+2v'u+\,v''}{ \left( {u}^{2}+1 \right)u}}\\
\end{equation}
that admits the three-dimensional point symmetry algebra generated
by
\begin{equation}\label{12e2}
\begin{array}{lll}
Y_1=\frac{\partial}{\partial v},&Y_2=\frac{1}{u}\frac{\partial}{\partial v},&Y_3=(u^2+1)\frac{\partial}{\partial u}+\frac{1}{u}(4uv-v-\sqrt{u^2+1})\frac{\partial}{\partial v},\\
\end{array}
\end{equation}
with the nonzero commutators
\begin{equation}\label{12e3}
\begin{tabular}{llllll}
$[Y_{1}, Y_{3}] =4 Y_{1}-Y_{2}$,& $[Y_{2}, Y_{3}] =Y_{1}+4Y_{2}$.\\
\end{tabular}
\end{equation}
Here dim $G'=2$,  rank $(G')=1$ and the adjoint action of
$G/G'=<\overline{Y}_3>$ on  $G'=<Y_1,Y_2>$ is given by
\begin{equation}\label{12e4}
{\rm ad}\,(\overline{Y}_3)= \left({
\begin{array}{cc}
-4 & -1 \\
1  &-4\\
\end{array}
}\right).
\end{equation}
The eigenvalues are $-4\pm \sqrt{-1}$ with eigenvectors $Y_2\pm
\sqrt{-1} ~Y_1$ respectively. Using the Theorem \ref{th1}, the
third-order ODE (\ref{12e1}) can be transformed to the canonical
form of $L^{II}_{3:7}$ via a point transformation.

To construct such a point transformation, one needs to match the
real part and imaginary part of an eigenvector of ad $(G/G')$ on
$G'$ with $X=\frac{\partial }{\partial y}$ and $Y=x\frac{\partial
}{\partial y}$ in the following way:
\begin{equation}\label{12e5}
\begin{array}{cc}
X=Y_2,&Y=Y_1.\\
\end{array}
\end{equation}
Applying the correspondence (\ref{12e5}) to the point
transformation $u=\phi(x,y), v=\psi(x,y)$ yields the system
\begin{equation}\label{12e6}
\begin{array}{cccc}
\phi_y=0,&x\phi_y=0,&\psi_y=\frac{1}{\phi},&x\psi_y=1.\\
\end{array}
\end{equation}
Solving the system (\ref{12e6}) gives the following point
transformation
\begin{equation}\label{12e7}
\begin{array}{cc}
u=x,& v=\frac{y}{x},\\
\end{array}
\end{equation}
which transforms the vector $Y_3$ which is linearly independent of
$G'$ to
\begin{equation}\label{12e8}
\begin{array}{cc}
(\frac{1}{\sin \theta})Z=(1+x^2)\, \partial_{x}+(y(x-b)+f(x))\,\partial_{y},\\
\end{array}
\end{equation}
with $f(x)=-\sqrt{x^2+1}$ and $b=\cot{\theta}=-4$. The function
$f(x)$ can be absorbed using the transformation
\begin{equation}\label{12e9}
\begin{array}{lll}
\tilde{x}=x,& \tilde{y}&=y-\sqrt{x^2+1} ~e^{-b\tan^{-1}x}\int \frac{e^{b \tan^{-1}x}}{{(x^2+1)}^{\frac{3}{2}}} f(x) dx=y-\frac{1}{4}\sqrt{x^2+1}.\\
\end{array}
\end{equation}
Finally, the composition of the transformations (\ref{12e7}) and
(\ref{12e9}) transforms ODE (\ref{12e1}) after simplification
using the identity $\tan^{-1}
x+\tan^{-1}\frac{1}{x}=\frac{\pi}{2}$ to its canonical form
\begin{equation}\label{12e10}
\tilde{y}'''=\frac{\tilde{y}''}{1+\tilde{x}^2}\left(f\left(\tilde{y}''{(\tilde{x}^2+1)}^{\frac{3}{2}}e^{b\tan^{-1}\tilde{x}}\right)-3\tilde{x}\right)\\
\end{equation}
 with $f(z)=-\frac{e^{4\pi}}{z^2}$ and $b=-4$.
\end{example}
%%%%%%%%%%%%%%%%%%%%%%%%%%%%%%%%%%%%%%%%%%%%%%%%%%%%%%%%%%%%%%%%%%%%%%%%%%%%%%%%%%%%%%%%%%%%%%%%%%%%%%%%%%%%%%%%%%%%%%%%%%%%%%
\begin{example}\rm {$L^{I}_{3:8}$}\label{ex13}\\
 Consider the ODE
\begin{equation}\label{13e1}
v''= {(v+v')}^{3}-\frac{1}{2}v-\frac{3}{2}v'\\
\end{equation}
that admits the three-dimensional point symmetry algebra generated
by
\begin{equation}\label{13e2}
\begin{array}{lll}
Y_1=\frac{\partial}{\partial u},&Y_2=\exp(-u)\frac{\partial}{\partial v},&Y_3=v\exp(u)\frac{\partial}{\partial u}-\frac{1}{2}v^2\exp(u)\frac{\partial}{\partial v},\\
\end{array}
\end{equation}
with the nonzero commutators
\begin{equation}\label{13e3}
\begin{tabular}{llllll}
$[Y_{1}, Y_{2}] = -Y_{2}$,& $[Y_{1}, Y_{3}] = Y_{3}$,& $[Y_{2}, Y_{3}] = Y_{1}$.\\
\end{tabular}
\end{equation}
Since dim $G'=3$, the Killing form is indefinite and rank $G=2$,
using the Theorem \ref{th1}, the second-order ODE (\ref{13e1}) can
be transformed to one of  the three canonical forms $L^{I}_{3:8}$,
$L^{II}_{3:8}$ and $L^{III}_{3:8}$ via a point transformation.

Since the eigenvalues of ad $(Y_1)$ are  $\pm 1$, then by a
scaling, as explained in section \ref{4.3}, one can get the change
of basis
\begin{equation}\label{13e4}
\begin{array}{cccc}
X=Y_3,& Y= -2Y_2,&Z=2Y_1.\\
\end{array}
\end{equation}
This maps the nonzero commutators (\ref{13e3}) to the standard
relations given by $$[Z,X]=2X,[Z,Y]=-2Y,[X,Y]=Z.$$ Applying the
correspondence (\ref{13e4}) to the point transformation
$u=\phi(x,y), v=\psi(x,y)$ yields the system
\begin{equation}\label{13e5}
\begin{array}{llll}
\phi_y=\exp(\phi)\psi,&\psi_y=-\frac{1}{2}\exp(\phi){\psi}^2,\\
-2xy\phi_x+(\epsilon~x^2-y^2)\phi_y=0,&-2xy\psi_x+(\epsilon~x^2-y^2)\psi_y=-2\exp(-\phi),\\
x\phi_x+y\phi_y=-1,&x\psi_x+y\psi_y=0,\\
\end{array}
\end{equation}
for some $\epsilon \in \{0,1,-1\}$.
Since the system (\ref{13e5}) is consistent for $\epsilon=0$, its
solution
\begin{equation}\label{13e6}
\begin{array}{cc}
u=\ln{\left(\frac{x}{y^2}\right)}+c_1,& v=-2e^{-c_1}\left(\frac{y}{x}\right),\\
\end{array}
\end{equation}
transforms ODE (\ref{13e1}) to the canonical form $L^{I}_{3:8}$
\begin{equation}\label{13e7}
xy''=C{y'}^3-\frac{1}{2}~y',\\
\end{equation}
with $C=4e^{-2c_1}$.
\end{example}
%%%%%%%%%%%%%%%%%%%%%%%%%%%%%%%%%%%%%%%%%%%%%%%%%%%%%%%%%%%%%%%%%%%%%%%%%%%%%%%%%%%%%%%%%%%%%%%%%%%%%%%%%%%%%%%%%%%%%%%%%%%%%%
\begin{example}\rm {$L^{II}_{3:8}$} \label{ex14}\\
 Consider the ODE
\begin{equation}\label{14e1}
u^4 v v''= {({v'}^2 u^4+1)}^{\frac{3}{2}}-1-2v v' u^3-{v'}^2 u^4\\
\end{equation}
that admits the three-dimensional point symmetry algebra generated
by
\begin{equation}\label{14e2}
\begin{array}{lll}
Y_1=u^2\frac{\partial}{\partial u},&Y_2=-u\frac{\partial}{\partial u}+v\frac{\partial}{\partial v},&Y_3=(u^2v^2-1)\frac{\partial}{\partial u}+\frac{2v}{u}\frac{\partial}{\partial v},\\
\end{array}
\end{equation}
with the nonzero commutators
\begin{equation}\label{14e3}
\begin{tabular}{llllll}
$[Y_{1}, Y_{2}] = Y_{1}$,& $[Y_{1}, Y_{3}] = -2Y_{2}$,& $[Y_{2}, Y_{3}] = Y_{3}$.\\
\end{tabular}
\end{equation}
Since dim $G'=3$, the Killing form is indefinite and rank $G=2$,
using the Theorem \ref{th1}, the second-order ODE (\ref{14e1}) can
be transformed to one of  the three canonical forms $L^{I}_{3:8}$,
$L^{II}_{3:8}$ and $L^{III}_{3:8}$ via a point transformation.

Since the eigenvalues of ad $(Y_2)$ are  $\pm 1$, by a scaling, as
explained in section \ref{4.3}, one can get the change of basis
\begin{equation}\label{14e4}
\begin{array}{cccc}
X=Y_3,& Y= Y_1,&Z=2Y_2.\\
\end{array}
\end{equation}
This maps the nonzero commutators (\ref{14e3}) to the standard
relations given by $$[Z,X]=2X,[Z,Y]=-2Y,[X,Y]=Z.$$ Applying the
correspondence (\ref{14e4}) to the point transformation
$u=\phi(x,y), v=\psi(x,y)$ yields the system
\begin{equation}\label{14e5}
\begin{array}{llll}
\phi_y={\phi}^2{\psi}^2-1,&\psi_y=\frac{2\psi}{\phi},\\
-2xy\phi_x+(\epsilon~x^2-y^2)\phi_y={\phi}^2,&-2xy\psi_x+(\epsilon~x^2-y^2)\psi_y=0,\\
x\phi_x+y\phi_y=\phi,&x\psi_x+y\psi_y=-\psi,\\
\end{array}
\end{equation}
for some $\epsilon \in \{0,1,-1\}$.
Since the system (\ref{14e5}) is consistent for $\epsilon=1$, its
solution
\begin{equation}\label{14e6}
\begin{array}{cc}
u=-\frac{x^2+y^2}{y},& v=\frac{x}{x^2+y^2},\\
\end{array}
\end{equation}
transforms ODE (\ref{14e1}) to the canonical form $L^{II}_{3:8}$
\begin{equation}\label{14e7}
xy''=y'+{y'}^3+C(1+{y'}^2)^{\frac{3}{2}},\\
\end{equation}
with $C=1$.
\end{example}
%%%%%%%%%%%%%%%%%%%%%%%%%%%%%%%%%%%%%%%%%%%%%%%%%%%%%%%%%%%%%%%%%%%%%%%%%%%%%%%%%%%%%%%%%%%%%%%%%%%%%%%%%%%%%%%%%%%%%%%%%%%%%%
\begin{example}\rm {$L^{III}_{3:8}$} \label{ex15}\\
 Consider the ODE
\begin{equation}\label{15e1}
u^4 v^3 v''= 3v^2{v'}^2u^4-v^6-2v^3v'u^3+{({v'}^2u^4-v^4)}^{\frac{3}{2}}\\
\end{equation}
that admits the three-dimensional point symmetry algebra generated
by
\begin{equation}\label{15e2}
\begin{array}{lll}
Y_1=u^2\frac{\partial}{\partial u},&Y_2=-u\frac{\partial}{\partial u}-v\frac{\partial}{\partial v},&Y_3=(1+\frac{u^2}{v^2})\frac{\partial}{\partial u}+\frac{2v}{u}\frac{\partial}{\partial v},\\
\end{array}
\end{equation}
with the nonzero commutators
\begin{equation}\label{15e3}
\begin{tabular}{llllll}
$[Y_{1}, Y_{2}] = Y_{1}$,& $[Y_{1}, Y_{3}] = 2Y_{2}$,& $[Y_{2}, Y_{3}] = Y_{3}$.\\
\end{tabular}
\end{equation}
Here again the Killing form is non-degenerate and indefinite and
rank $G=2$. Using the Theorem \ref{th1}, the second-order ODE
(\ref{15e1}) can be transformed to one of  the three canonical
forms $L^{I}_{3:8}$, $L^{II}_{3:8}$ and $L^{III}_{3:8}$ via a
point transformation.

Since the eigenvalues of ad $(Y_2)$ are  $\pm 1$, then by a
scaling, as explained in section \ref{4.3}, one has the change of
basis
\begin{equation}\label{15e4}
\begin{array}{cccc}
X=Y_3,& Y= -Y_1,&Z=2Y_2.\\
\end{array}
\end{equation}
This maps the nonzero commutators (\ref{15e3}) to the standard
relations given by $$[Z,X]=2X,[Z,Y]=-2Y,[X,Y]=Z.$$ Applying the
correspondence (\ref{15e4}) to the point transformation
$u=\phi(x,y), v=\psi(x,y)$ yields the system
\begin{equation}\label{15e5}
\begin{array}{llll}
\phi_y=1+\frac{{\phi}^2}{{\psi}^2},&\psi_y=\frac{2\psi}{\phi},\\
-2xy\phi_x+(\epsilon~x^2-y^2)\phi_y=-{\phi}^2,&-2xy\psi_x+(\epsilon~x^2-y^2)\psi_y=0,\\
x\phi_x+y\phi_y=\phi,&x\psi_x+y\psi_y=\psi,\\
\end{array}
\end{equation}
for some $\epsilon \in \{0,1,-1\}$.
Since the system (\ref{15e5}) is consistent for $\epsilon=-1$, its
solution
\begin{equation}\label{15e6}
\begin{array}{cc}
u=\frac{y^2-x^2}{y},& v=\frac{y^2-x^2}{x},\\
\end{array}
\end{equation}
transforms ODE (\ref{15e1}) to the canonical form $L^{III}_{3:8}$
\begin{equation}\label{15e7}
xy''=y'-{y'}^3+C(1-{y'}^2)^{\frac{3}{2}},\\
\end{equation}
with $C=-1$.
\end{example}
%%%%%%%%%%%%%%%%%%%%%%%%%%%%%%%%%%%%%%%%%%%%%%%%%%%%%%%%%%%%%%%%%%%%%%%%%%%%%%%%%%%%%%%%%%%%%%%%%%%%%%%%%%%%%%%%%%%%%%%%%%%%%%
\begin{example}\rm {$L^{IV}_{3:8}$} \label{ex16}\\
 Consider the ODE
\begin{equation}\label{16e1}
v'''=\frac{3}{2}\frac{{v''}^2}{v'}-\frac{{v'}^3}{v^2}\\
\end{equation}
that admits the three-dimensional point symmetry
$\mathbf{subalgebra}$ generated by
\begin{equation}\label{16e2}
\begin{array}{lll}
Y_1=u\frac{\partial}{\partial u},&Y_2=\frac{\partial}{\partial u},&Y_3=\frac{1}{2}u^2\frac{\partial}{\partial u},\\
\end{array}
\end{equation}
with the nonzero commutators
\begin{equation}\label{16e3}
\begin{tabular}{llllll}
$[Y_{1}, Y_{2}] = -Y_{2}$,& $[Y_{1}, Y_{3}] = Y_{3}$,& $[Y_{2}, Y_{3}] = Y_{1}$.\\
\end{tabular}
\end{equation}
Here dim $G'=3$, the Killing form is indefinite and rank $G=1$.
Using the Theorem \ref{th1}, the third-order ODE (\ref{16e1}) can
be transformed to the canonical form $L^{IV}_{3:8}$ via a point
transformation.

Since the eigenvalues of ad $(Y_1)$ are $\pm 1$, then by a scaling,
as explained in section \ref{4.3}, one can get the change of basis
\begin{equation}\label{16e4}
\begin{array}{cccc}
X=Y_3,& Y= -2Y_2,&Z=2Y_1.\\
\end{array}
\end{equation}
This maps the nonzero commutators (\ref{13e3}) to the standard
relations given by $$[Z,X]=2X,[Z,Y]=-2Y,[X,Y]=Z.$$ Applying the
correspondence (\ref{16e4}) to the point transformation
$u=\phi(x,y), v=\psi(x,y)$ yields the system
\begin{equation}\label{16e5}
\begin{array}{llll}
\phi_y=\frac{1}{2}{\phi}^2,&\psi_y=0,\\
-y^2\phi_y=-2,&-y^2\psi_y=0,\\
y\phi_y=-\phi,&y\psi_y=0.\\
\end{array}
\end{equation}
The solution of the system (\ref{16e5}) gives a point
transformation
\begin{equation}\label{16e6}
\begin{array}{cc}
u=-\frac{2}{y},& v=x,\\
\end{array}
\end{equation}
that transforms ODE (\ref{16e1}) to the canonical form
$L^{IV}_{3:8}$
\begin{equation}\label{16e7}
y'''=\frac{3}{2}\frac{{y''}^2}{y'}+f(x)y',\\
\end{equation}
with $f(x)=\frac{1}{x^2}$.
\end{example}
%%%%%%%%%%%%%%%%%%%%%%%%%%%%%%%%%%%%%%%%%%%%%%%%%%%%%%%%%%%%%%%%%%%%%%%%%%%%%%%%%%%%%%%%%%%%%%%%%%%%%%%%%%%%%%%%%%%%%%%%%%%%%%
\begin{example}\rm {$L_{3:9}$} \label{ex17}\\
 Consider the ODE
\begin{equation}\label{17e1}
v''=-{v'}^3\cos u \sin u-2v'\cot u+\csc u{({v'}^2\sin^2 u +1)}^\frac{3}{2}\\
\end{equation}
that admits the three-dimensional point symmetry algebra generated
by
\begin{equation}\label{17e2}
\begin{array}{lll}
Y_1=\frac{\partial}{\partial v},&Y_2=\sin v\frac{\partial}{\partial u}+\cos v\cot u\frac{\partial}{\partial v},&Y_3=\cos v\frac{\partial}{\partial u}-\sin v\cot u\frac{\partial}{\partial v},\\
\end{array}
\end{equation}
with the nonzero commutators
\begin{equation}\label{17e3}
\begin{tabular}{llllll}
$[Y_{1}, Y_{2}] =Y_{3}$,& $[Y_{1}, Y_{3}] =-Y_{2}$,& $[Y_{2}, Y_{3}] =Y_{1}$.\\
\end{tabular}
\end{equation}
Since the Killing form is negative definite, using Theorem
\ref{th1}, the second-order ODE (\ref{17e1}) can be transformed to
the canonical form of $L_{3:9}$ via a point transformation.

In order to construct such a point transformation, pick any vector
like $Y_1$ and find its non-zero eigenvalues. Here ad $(Y_1)$ has
$\pm \sqrt{-1}$ as  eigenvalues with eigenvectors $Y_3\pm
\sqrt{-1}~ Y_2$ respectively. One needs to match the vector $Y_1$
with $X$ and a multiple of eigenvector $Y_3+ \sqrt{-1}~ Y_2$ with
the vector $Y-\sqrt{-1}~Z$ such that $[Y,Z]=X$ in the following
way:
\begin{equation}\label{17e5}
\begin{array}{ccc}
X=Y_1,&Y=Y_3,&Z=-Y_2.\\
\end{array}
\end{equation}
Applying the correspondence (\ref{17e5}) to the point
transformation $u=\phi(x,y), v=\psi(x,y)$ yields the system
\begin{equation}\label{17e6}
\begin{array}{llll}
\phi_x=0,&\psi_x=1,\\
y\sin x ~\phi_x  + \left( {y}^{2}+1\right) \cos x ~\phi_y =\cos \psi,&y\sin x ~\psi_x + \left( {y}^{2}+1 \right) \cos x ~\psi_y =-\sin  \psi  \cot \phi,\\
y\cos x ~\phi_x  - \left( {y}^{2}+1\right) \sin x ~\phi_y  =-\sin\psi,&y\cos x ~\psi_x - \left( {y}^{2}+1 \right) \sin x ~\psi_y  =-\cos \psi \cot \phi.\\
\end{array}
\end{equation}
Solution of the system (\ref{17e6}) gives the required point
transformation. A systematic way of solving such a nonlinear
system is as follows:

One can match the vector $Y_1$ with $X$ through the canonical
coordinates of $Y_1$ as
\begin{equation}\label{17e7}
\begin{array}{cc}
u=y,& v=x.\\
\end{array}
\end{equation}
This transforms the vector $Y-\sqrt{-1}~Z$ to
\begin{equation}\label{17e8}
\begin{array}{cc}
Y-\sqrt{-1}~Z=Y_3+\sqrt{-1}~Y_2=e^{\sqrt{-1}~x}\left[ \left( f_1(y)\partial_{x}+f_2(y)\partial_{y}\right) + \sqrt{-1} \left( f_3(y)\partial_{x}+f_4(y)\partial_{y}\right)\right],\\
\end{array}
\end{equation}
with $f_1(y)=0,f_2(y)=1,f_3(y)=\cot{y}$ and $f_4(y)=0$. Now using
the formula (\ref{te4}), the vector $Y-\sqrt{-1}~Z$ can be
transformed using the transformation
\begin{equation}\label{17e9}
\begin{array}{lll}
\tilde{x}=x+{\tan}^{-1}\left(\frac{f_4}{f_2}\right)=x,& \tilde{y}&=\frac{f_1f_4-f_2f_3}{\sqrt{f_2^2+f_4^2}}=-\cot y.\\
\end{array}
\end{equation}
to the canonical form
\begin{equation}\label{17e10}
\begin{array}{cc}
Y-\sqrt{-1}~Z=e^{\sqrt{-1}\tilde{x}}\left[  (1+{\tilde{y}}^2)\partial_{\tilde{y}} - \sqrt{-1}~\tilde{y}\partial_{\tilde{x}}\right],\\
\end{array}
\end{equation}
Finally, the composition of the transformations (\ref{17e7}) and
(\ref{17e9}) given by
\begin{equation}\label{17e11}
\begin{array}{cc}
u=-\cot^{-1}{\tilde{y}},& v=\tilde{x},\\
\end{array}
\end{equation}
transforms ODE (\ref{17e1}) to its canonical form
\begin{equation}\label{17e12}
\tilde{y}''=C{\left( \frac{\tilde{y}^{'2}+\tilde{y}^2+1}{1+\tilde{y}^2}\right)}^{\frac{3}{2}}-\tilde{y}\\
\end{equation}
 with $C=1$.

Moreover, a solution of the nonlinear system (\ref{17e6}) is
\begin{equation}\label{17e13}
\begin{array}{cc}
u=-\cot^{-1}{y},& v=x.\\
\end{array}
\end{equation}
\end{example}
%%%%%%%%%%%%%%%%%%%%%%%%%%%%%%%%%%%%%%%%%%%%%%%%%%%%%%%%%%%%%%%%%%%%%%%%%%%%%%%%%%%%%%%%%%%%%%%%%%%%%%%%%%%%%%%%%%%%%%%%%%%%%%
\section{Conclusion}

The Lie-Bianchi classification of three-dimensional algebras and their
realizations as vector fields in ${\mathbb{R}}^2$ are recovered in an
algorithmic way. This is done in such a way that one can read off the
type of the algebra from its invariants like the dimension of its
commutator or the centralizer of its commutator and its rank.

The compact and non-compact Lie algebras are treated uniformly in
a manner that makes their realizations as vector fields in the
plane transparent.

The algorithms are illustrated by examples for each type of three
dimensional algebras. The procedures works in principle for any ODE
which admits a three-dimensional subalgebra of symmetries.
%-----------------------------------------------------------------------------%
\subsection*{Acknowledgments}
%-----------------------------------------------------------------------------%
 The authors would like to thank the King Fahd University of
Petroleum and Minerals for its support and excellent research
facilities. FM is grateful to the NRF of South Africa for research funding support.
%-----------------------------------------------------------------------------%

\end{document}